\documentclass[bj,preprint]{imsart}
\usepackage{amssymb}
\usepackage{epsfig}
\usepackage[french,english]{babel}
\usepackage{graphicx}
\usepackage[latin1]{inputenc}
{\catcode `\@=11 \global\let\AddToReset=\@addtoreset}
\AddToReset{equation}{section}

\newtheorem{theorem}{Theorem}[section]
\newtheorem{lemma}{\bf Lemma}[section]
\newtheorem{corollary}{\bf Corollary}[section]
\newtheorem{property}{Property}[section]
\newtheorem{proposition}{Proposition}[section]
\newtheorem{@definition}{\sc Definition}[section]

\newtheorem{@remark}{\sc Remark}[section]
\newenvironment{remark}{\begin{@remark}\rm}{\end{@remark}}
\newtheorem{@example}{\sc Example}[section]
\newenvironment{example}{\begin{@example}\rm}{\end{@example}}

\newcommand{\beqn}{\begin{displaymath}}
\newcommand{\eeqn}{\end{displaymath}}
\newcommand{\beq}{\begin{equation}}  
\newcommand{\eeq}{\end{equation}}

\def\mathsf{\bf}
\def\N{\mathbb{N}}

\def\L{\mathbb{L}}
\def\R{\mathbb{R}}
\def\Z{\mathbb{Z}}

\def\i{\mathrm i}
\def\d{\mathrm d}
\def\e{\mathrm e}

\def\E{\mathrm E}
\def\P{\mathrm P}

\def\text{\mbox}

\def\1{{\bf 1}}

\newcommand{\mbf}[1]{\mbox{\boldmath $#1$}}

\newcommand{\Var}{\mbox{var}}
\newcommand{\Cov}{\mbox{cov}}

\def\simn{\renewcommand{\arraystretch}{0.5}
\begin{array}[t]{c}
\stackrel{}{\sim} \\
{\scriptstyle n\rightarrow\infty}
\end{array}\renewcommand{\arraystretch}{1}}

\def\limiteloin{\renewcommand{\arraystretch}{0.5}
\begin{array}[t]{c}
\stackrel{{\cal D}}{\longrightarrow} \\
{\scriptstyle n \rightarrow\infty}
\end{array}\renewcommand{\arraystretch}{1}}

\def\limiteloit{\renewcommand{\arraystretch}{0.5}
\begin{array}[t]{c}
\stackrel{{\cal D}}{\longrightarrow} \\
{\scriptstyle t \rightarrow 0}
\end{array}\renewcommand{\arraystretch}{1}}

\def\limiteproban{\renewcommand{\arraystretch}{0.5}
\begin{array}[t]{c}
\stackrel{{\cal P}}{\longrightarrow} \\
{\scriptstyle n \rightarrow\infty}
\end{array}\renewcommand{\arraystretch}{1}}

\def\limiten{\renewcommand{\arraystretch}{0.5}
\begin{array}[t]{c}
\stackrel{}{\longrightarrow} \\
{\scriptstyle n\rightarrow\infty}
\end{array}\renewcommand{\arraystretch}{1}}

\def\limitepsn{\renewcommand{\arraystretch}{0.5}
\begin{array}[t]{c}
\stackrel{a.s.}{\longrightarrow} \\
{\scriptstyle n\rightarrow\infty}
\end{array}\renewcommand{\arraystretch}{1}}

\def\limitefddl{\renewcommand{\arraystretch}{0.5}
\begin{array}[t]{c}
\stackrel{f.d.d.}{\longrightarrow} \\
{\scriptstyle \delta \rightarrow 0}
\end{array}\renewcommand{\arraystretch}{1}}

\def\limitefdd{\renewcommand{\arraystretch}{0.5}
\begin{array}[t]{c}
\stackrel{f.d.d.}{\longrightarrow} \\
\end{array}\renewcommand{\arraystretch}{1}}

\def\limitefddn{\renewcommand{\arraystretch}{0.5}
\begin{array}[t]{c}
\stackrel{f.d.d.}{\longrightarrow} \\
{\scriptstyle n\rightarrow\infty}
\end{array}\renewcommand{\arraystretch}{1}}

\def\limitet{\renewcommand{\arraystretch}{0.5}
\begin{array}[t]{c}
\stackrel{}{\longrightarrow} \\
{\scriptstyle t\rightarrow\infty}
\end{array}\renewcommand{\arraystretch}{1}}

\def\limitet0{\renewcommand{\arraystretch}{0.5}
\begin{array}[t]{c}
\stackrel{}{\longrightarrow} \\
{\scriptstyle t\rightarrow 0}
\end{array}\renewcommand{\arraystretch}{1}}

\newtheorem{thm}{Theorem}

\newtheorem{lem}{Lemma}

\begin{document}
\begin{frontmatter}

\title{Measuring the roughness of random paths by increment ratios}
\runtitle{Measuring the roughness of random paths}
\begin{aug}
\author{Jean-Marc Bardet\thanksref{t1} \ead[label=e1]{bardet@univ-paris1.fr}
  }
  \and
  \author{Donatas~Surgailis \thanksref{t2} \ead[label=e2]{sdonatas@ktl.mii.lt}
  \ead[label=u2,url]{http://www.mii.lt/index.php}}
\runauthor{J.-M. Bardet and D. Surgailis}

\affiliation{University Panth\'eon-Sorbonne, Institute of Mathematics and Informatics}

\address{SAMM,
University Paris 1,
90 Rue de Tolbiac\\
75634 Paris Cedex 13, FRANCE\\
\printead{e1}}

\address{Institute of Mathematics and Informatics\\ Akademijos 4, 08663 Vilnius, LITHUANIA \\ \printead{e2}}
\thankstext{t1}{This author thanks the Steklov Institute of St Petersburg (Russia) for a nice and fruitful invitation}
\thankstext{t2}{Research of this author was  supported by the Lithuanian State Science and Studies
Foundation grant T-70/09.}
\end{aug}
\begin{abstract}
~~A statistic based on increment ratios (IR) and related to zero  crossings of increment sequence is defined and studied for
measuring the roughness of random paths. The main advantages of this statistic are 
robustness to smooth additive and multiplicative trends and applicability
to infinite variance processes. The existence of the IR statistic limit (called the IR-roughness below) is closely related
to the existence of a tangent process. Three particular cases where the IR-roughness exists and is
explicitly computed are considered. Firstly, for a diffusion
process with smooth diffusion and drift
coefficients, the IR-roughness coincides with the IR-roughness of a Brownian motion and its convergence rate is obtained.
Secondly, the case of rough
Gaussian processes is studied in detail under general
assumptions which do not require stationarity conditions.
Thirdly, the IR-roughness of a Lévy process with $\alpha-$stable tangent process is established and
can be used to estimate the fractional parameter $\alpha \in (0,2)$ following a central limit theorem.
\end{abstract}
\begin{keyword}[class=AMS]
\kwd[Primary ]{62H12, 60G17;} \ \kwd[Secondary ]{62F12, 60G18, 60G15, 60G51}
\end{keyword}

\begin{keyword}
\kwd{Estimation of the local regularity function of
stochastic process; Limit theorems; Hölder exponent; Tangent process; Zero crossings; Diffusion processes; Fractional Brownian motion; Multifractional Brownian motion; Lévy processes}
\end{keyword}

\end{frontmatter}

\section{Introduction and the main results}

It is well-known that random functions are typically ``rough''
(non-differentiable), which raise the question of determining and
measuring roughness. Probably, the most studied roughness measures
are the Hausdorff dimension and the $p-$variation index. There exists
a considerable literature on statistical estimation of these and
related quantities from a discrete grid. Hence, different estimators
of the Hausdorff dimension have been studied, as the box-counting
estimator (see Hall and Wood, 1993 for stationary Gaussian processes
or L\'evy-V\'ehel and Peltier, 1994, for Gaussian processes with
stationary increments). To our knowledge, the $H$-variation
estimator, where $H$ is a measurable function, was first
proposed by Guyon and Leon (1989) for stationary Gaussian processes
where central and non-central limit theorems are established
following the Hermite rank of $H$ and the asymptotic local
properties of the variogram and its second derivative. Further
studies provided a continuation of this seminal paper in different
ways. Istas and Lang~(1997) studied  generalized
quadratic variations of Gaussian processes with stationary
increments. Coeurjolly~(2001 and 2005) studied
$\ell^p$-variations of fractional Brownian
motion and $\ell^2$-variations of multifractional
Brownian motion.
Coeurjolly~(2007) discussed $L$-variations
based on linear combinations of empirical quantiles for Gaussian
locally self-similar processes. An estimator counting the number
level crossings was investigated by Feuerverger {\it et al.}~(1994) for stationary Gaussian processes.

\smallskip

In the present paper we introduce a new characteristic of roughness,
defined as a sum of ratios of consecutive increments. For a
real-valued function $f = (f(t), t \in [0,1]) $, define recursively
\begin{eqnarray}
\Delta^{1,n}_jf&:=&
f \big( \frac{j+1}{n}\big) -  f \big( \frac{j}{n}\big), \nonumber \\
\Delta^{p,n}_j f&:=&\Delta^{1,n}_j \Delta^{p-1,n}_j f
\  = \  \sum_{i=0}^p (-1)^{p-i} {p \choose i} f\big(\frac{j+i}{n}\big), \label{Delta-p}
\end{eqnarray}
so that $\Delta^{p,n}_j f$ denotes the  $p-$order increment of $f$ at $\frac{j}{n}$, $\, p=1,2, \cdots,  \ j=0,1, \cdots, n-p$.  Let
\begin{eqnarray}
R^{p,n} (f) &:=&\frac{1}{n-p} \sum_{k=0}^{n-p-1}
\frac{\big|\Delta^{p,n}_k f +  \Delta^{p,n}_{k+1} f \big|}
{|\Delta^{p,n}_k f| + |\Delta^{p,n}_{k+1} f| }, \label{R_p}
\end{eqnarray}
with the convention $\frac{0}{0}:=1$. In particular,
\begin{equation}
R^{1,n}(f) \ =  \  \frac{1}{n-1} \sum_{k=0}^{n-2}
\frac{\left|\mbox{$f\left(\frac{k+1}{n}\right) -
f\left(\frac{k}{n}\right) +f\left(\frac{k+2}{n}\right) -
f\left(\frac{k+1}{n}\right)$} \right|}
{\left|\mbox{$f\left(\frac{k+1}{n}\right) -
f\left(\frac{k}{n}\right)$}\right| +
\left|\mbox{$f\left(\frac{k+2}{n}\right) -
f\left(\frac{k+1}{n}\right)$} \right|}. \label{R_1}
\end{equation}
Note the ratio on the right-hand side of (\ref{R_p}) is either $1$ or
less than $1$ depending on whether the consecutive increments
$\Delta^{p,n}_k f$  and $\Delta^{p,n}_{k+1} f$ have same signs or
different signs; moreover, in the latter case, this ratio generally
is small  whenever the increments are similar in magnitude
(``cancel each other''). Clearly,  $0\le R^{p,n}(f) \le 1$ for any
$f, n, p$. Thus if $\lim R^{p,n}(f) $ exists when $n \to \infty$, the quantity
$R^{p,n}(f) $ can be used to estimate this limit which represents the  ``mean roughness of $f$'' also called the {\em $p-$th order IR-roughness of $f$} below. We show below that these definitions can be extended to sample paths of very general random processes, e.g. stationary
processes, processes with stationary and nonstationary increments,
and even $\L^q$-processes with $q<1$.

Let us describe the main results of this paper. Section \ref{general} derives some
general results on asymptotic behavior of this estimator.
Proposition~2.1 says that, for a sufficiently smooth function $f$,
the limit $\lim_{n \to \infty} R^{p,n}(f) =1$. In the most of the
paper, $f =X$ is a random process. Following Dobrushin~(1980), we
say that $X = (X_t,\, t \in {\R}) $ has a {\it small scale limit
$Y^{(t_0)}$ at point $t_0 \in {\R}$} if there exist a normalization
$A^{(t_0)}(\delta) \to \infty $ when $\delta \to 0$ and a random
process $Y^{(t_0)} = (Y^{(t_0)}_\tau, \tau \ge 0)$ such that
\begin{eqnarray}\label{dobrushin}
A^{(t_0)}(\delta) \left(X_{t_0+\tau \delta}- X_{t_0}\right) \limitefddl Y^{(t_0)}_\tau,
\end{eqnarray}
where $\limitefdd$ stands for weak convergence of finite dimensional
distributions. A related definition is given in Falconer~(2002,
2003) who called the limit process  $Y^{(t_0)}$ a {\it tangent
process} (at $t_0$). See also Benassi {\it et al.}~(1997). In many
cases, the normalization $A^{(t_0)} (\delta) = \delta^{H(t_0)} $,
where $0<H(t_0)< 1$ and the limit tangent process $Y^{(t_0)}$ is
self-similar with index $H(t_0)$ (Falconer, 2003 or Dobrushin,
1980). Proposition 2.2 states that if $X$ satisfies a similar
condition to (\ref{dobrushin}), then the statistic $R^{p,n} (X)$
converges to the integral
\begin{equation}
R^{p,n}(X) \limiteproban  \int_0^1 \E \left[ \frac{
|\Delta^p_0Y^{(t)} + \Delta^p_1 Y^{(t)}| }{|\Delta^p_0Y^{(t)}| +
|\Delta^p_1 Y^{(t)}|} \right]  \d t,  \label{la1}
\end{equation}
in probability, where $\Delta^p_j Y^{(t)} = \Delta^{p,1}_j Y^{(t)} =
\sum_{i=0}^p (-1)^{p-i} {p \choose i} Y^{(t)}_{j+i}, \  j=0,1$ is the corresponding increment of the tangent process
$Y^{(t)}$ at $t\in [0,1)$.  In the particular case when $X$ has
stationary increments, relation (\ref{la1}) becomes
\begin{equation}\label{la}
R^{p,n} (X) \limiteproban  \E \left[\frac{|\Delta^p_0 Y + \Delta^p_1
Y|}{|\Delta^p_0 Y| + |\Delta^p_1 Y|}\right].
\end{equation}
Section \ref{diffu} discusses the convergence in (\ref{la1}) for {\it diffusion
processes} $X$ admitting a stochastic differential $\d X = a_t \d
B(t) + b_t \d t$, where $B$ is a standard Brownian motion and
$(a_t), (b_t) $ are random (adapted) functions. It is clear that
under general regularity conditions on the diffusion and drift
coefficients $(a_t), (b_t) $, the process $X$ admits the same local
H\"older exponent as $B$ at each point $t_0 \in (0,1)$ and therefore
the IR-roughness of $X$ in (\ref{la1}) should not depend
on these coefficients and should coincide with the corresponding
limit for $X=B$. This is indeed the case since the tangent process
of $X$ at $t$ is easily seen to be $Y^{(t)} = a_{t} B$ and the
multiplicative factor $a_{t}$ cancels in the numerator and the
denominator of the fraction inside the expectation in (\ref{la1}).
See Proposition 3.1 for details, where the convergence rate
$O(n^{1/3})$ (a.s.) in (\ref{la1}) with explicit limit values
$\Lambda_p(1/2)$  is established for diffusions $X$ and $p=1,2$.

\smallskip

Considerable attention  is
given to the asymptotic behavior of the
statistic $R^{p,n}(X)$ for {\it ``fractal'' Gaussian
processes} (see Section \ref{gauss}). In such a frame, {\it fractional Brownian motion} (fBm
in the sequel) is a typical example. Indeed, if $X$ is a fBm with
parameter $H \in (0,1)$, then $X$ is also its self tangent process
for any $t \in [0,1]$ and (see Section \ref{gauss}):
\begin{eqnarray}\label{a.s.}
R^{p,n}(X) &\limitepsn& \Lambda_{p}(H) \qquad (p=1,2),
\end{eqnarray}
where
\begin{eqnarray}
\Lambda_p(H)&:=& \lambda (\rho_p(H)), \label{lambda_p}\\
\lambda(r)&:=&\frac{1}{\pi}\arccos (-r) +
\frac{1}{\pi}
\sqrt{\frac{1+r}{1-r}} \log \left(\frac{2}{1+r}\right), \label{lambda_0}\\
\displaystyle{\rho_p(H)}&:=&\displaystyle{{\rm corr}\big (\Delta^p_0
B_H, \Delta^p_1 B_H \big)},  \label{rho_p}
\end{eqnarray}
and where  $\Delta^1_j B_H = B_H(j+1)-B_H(j),~ \Delta^2_j B_H =
B_H(j+2)-2B_H(j+1)-B_H(j) \  (j\in {\Z})$ are respective increments
of fBm. Moreover,
\begin{equation}
\sqrt n \big (R^{p,n}(X)- \Lambda_{p}(H)\big )
\limiteloin {\cal N}(0,\Sigma_p(H)) \quad
\mbox{if}\ ~~~ \cases{p=1, \   \  0< H < 3/4, \cr
p=2, \   \    0<H<1, \cr}
\label{CLT_fBm}
\end{equation}
where $\limiteloin$ stands for weak convergence of probability distributions. The asymptotic variances $\Sigma_p(H)$ in (\ref{CLT_fBm}) are given by
\begin{equation}  \label{Sigma_p}
\Sigma_p(H):= \sum_{j\in {\Z}} \Cov\left(
\frac {\big | \Delta^p_0 B_{H}+\Delta^p_1B_{H}\big |}{\big |
\Delta^p_0 B_{H}\big |+\big |\Delta^p_1B_{H}\big |},  \frac {\big
| \Delta^p_j B_{H}+\Delta^p_{j+1}B_{H}\big |}{\big |
\Delta^p_j B_{H}\big |+\big |\Delta^p_{j+1} B_{H}\big |} \right).
\end{equation}
The graphs of $\Lambda_p(H)$ and $\sqrt{\Sigma_p(H)}\
(p=1,2) $ are given in Figures \ref{Figure1} and \ref{Figure3} below.
\begin{center}
\begin{figure}[h]
\begin{center}
\includegraphics[width=7 cm,height=5cm]{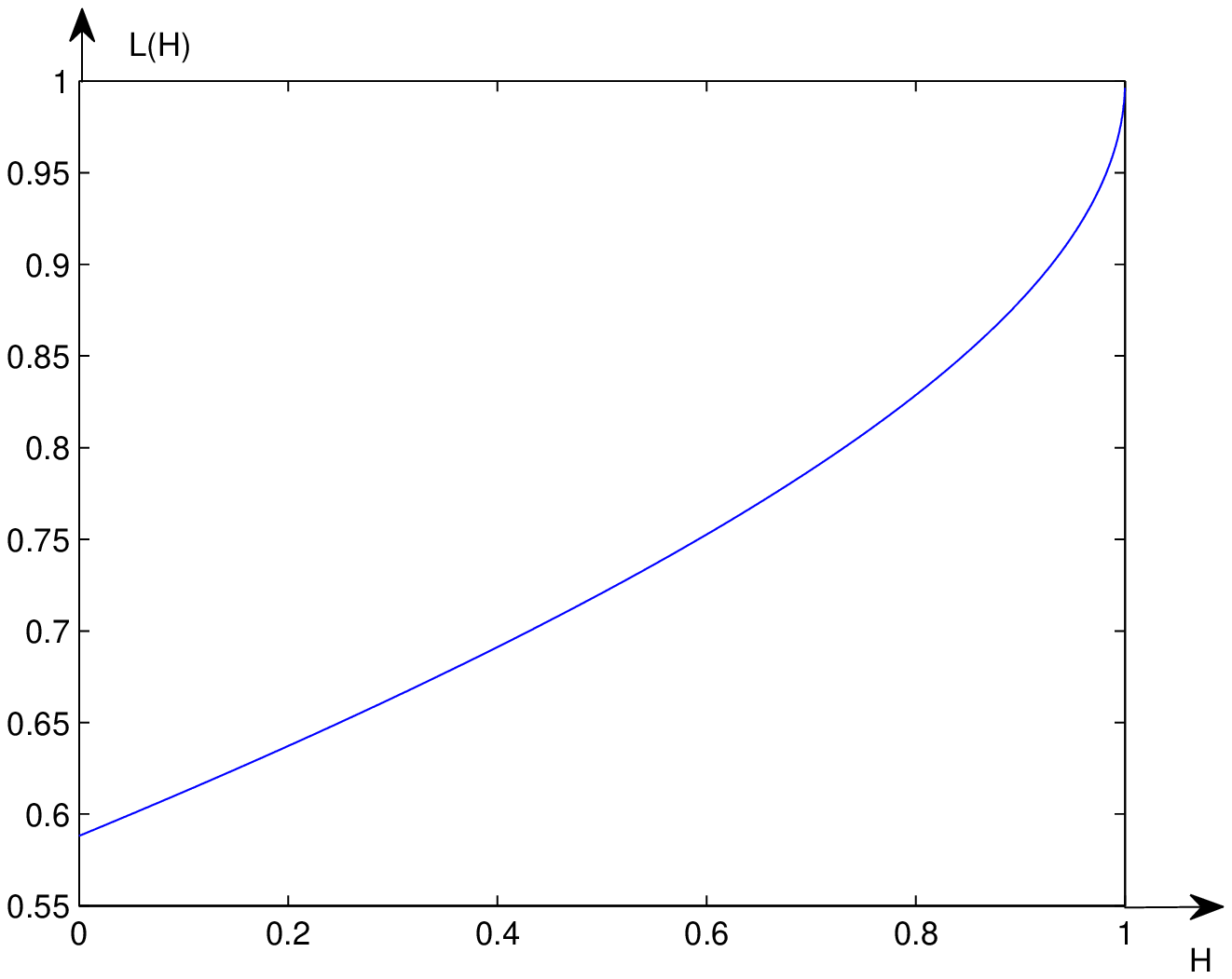}\hspace{2mm}
\includegraphics[width=7 cm,height=5cm]{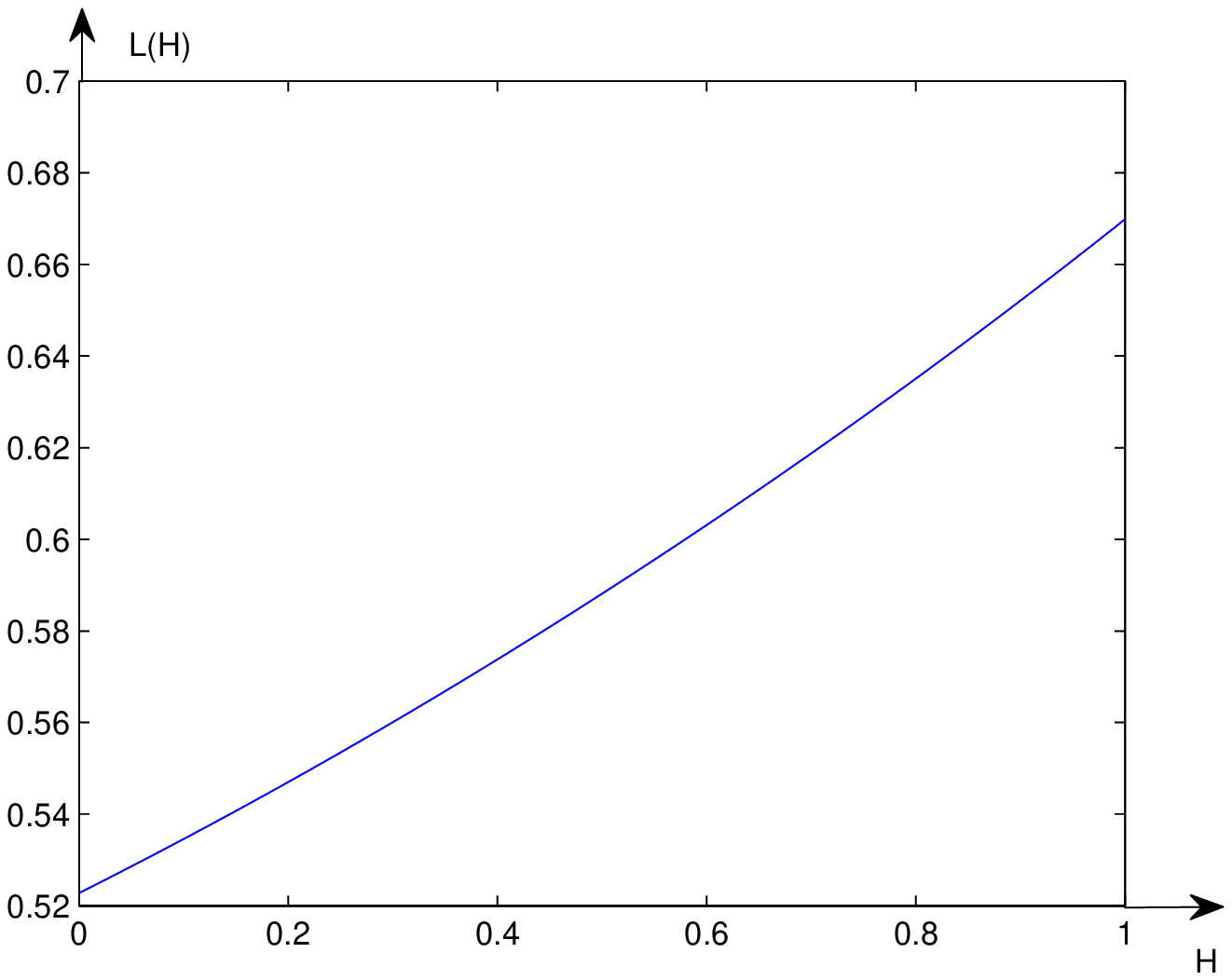}
\caption{The graphs of $\Lambda_1(H)$ (left) and
$\Lambda_2(H)$ (right). } \label{Figure1}
\end{center}
\end{figure}
\end{center}
\begin{center}
\begin{figure}[h]
\begin{center}
\includegraphics[width=10 cm,height=5cm]{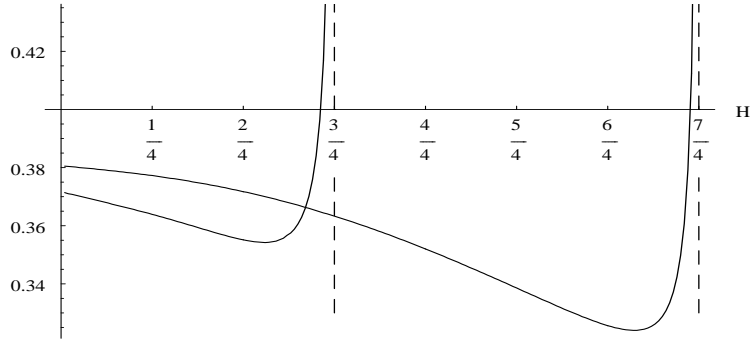}
\caption{The graphs of $\sqrt{\Sigma_p(H)}$, $p=1$ (with a pole at $3/4$) and $p=2$ (with a pole at $7/4$) {\rm (from
Stoncelis and Vai\v ciulis, 2008, with kind permission of the
authors)}} \label{Figure3}
\end{center}
\end{figure}
\end{center}
The difference in the range of the parameter $H$ for $p=1$ and $p=2$
in the central limit theorem in (\ref{CLT_fBm}) are due to the fact
that the second order increment process $\big(\Delta^2_j B_H, j\in
{\Z}\big)$ is a short memory stationary Gaussian process for any
$H \in (0,1)$, in contrast to the first order increment process
$\big(\Delta^1_j B_H, j \in {\Z}\big)$ which has long memory for $H>3/4$.\\
Generalizations of (\ref{a.s.}) and (\ref{CLT_fBm})  to Gaussian
processes having nonstationary increments are proposed in
Section \ref{gauss}. Roughly speaking, $R^{p,n}(X),\ p=1,2$
converge a.s.  and satisfy a central limit theorem, provided for any $t \in
[0,1]$ the process $X$ admits a fBm with parameter $H(t)$ as a
tangent process (more precise assumptions {\bf (A.1)}, {\bf
(A.1)}$^\prime$ and {\bf (A.2)}$_p$ are provided in Section
\ref{gauss}). In such frames, the limits in (\ref{a.s.}) are $\int_0^1
\Lambda_{p}(H(t))\d t$ instead
of $\Lambda_{p}(H)$  and the asymptotic
variances in (\ref{CLT_fBm})
also change. The case of Gaussian processes
with stationary increments is discussed in detail and the results are used to define a $\sqrt{n}-$consistent estimator of
$H$, under semiparametric assumptions on the
asymptotic behavior of the variogram or the spectral density.  Bardet and Surgailis~(2010) study a punctual estimator of $H(t_0)$ obtained
from a localization around $t_0 \in (0,1)$ of the statistic
$R^{2,n}(X)$.

The main advantages of estimators of the type (\ref{R_p}) involving
a scaling invariant function of increments 
seem to be the following. Firstly,
the estimator $R^{p,n}(X)$ essentially depends on local regularity
of the process $X$  and not on possible ``multiplicative and additive
factors'' such as diffusion and drift coefficients in Section \ref{diffu} or
smoothly multiplicative and additive trended Gaussian processes, see
Proposition \ref{trends} of Section \ref{gauss}. This property is important when
dealing with financial data involving heteroscedasticity and
volatility clustering. Such a robustness property (also  satisfied by the estimators
based on generalized quadratic variations of wavelet coefficients)  represents a clear advantage versus classical parametric Whittle or semi-parametric log-periodogram estimators.
Secondly, the estimators in (\ref{R_p}) are
bounded functionals and have finite
moments of any order.
Section~\ref{Levy} discusses jump L\'evy processes, with the L\'evy measure regularly varying of fractional index $\alpha\in (0,2)$ at the origin.
Using a modification of (\ref{R_p}), we define a $\sqrt{n}-$consistent estimator
of $\alpha $, together with a central limit theorem, in a very general semiparametric frame.
This result is new and interesting because there exist very few papers providing consistent estimators of $\alpha$ (to our knowledge, the only comparable results have been established in (Belomestny, 2010) and (Ait Sahalia and Jacod, 2009) in a financial and somewhat different context). Finally, in the Gaussian case, using the approximation formula provided in Remark \ref{Hlinear}, an estimator of $H$ based on $R^{2,n}(X)$ can be extremely simply computed:
$$
\widehat H^{(2)}_n\simeq \frac 1 {0.1468} \Big (\frac{1}{n-2} \sum_{k=0}^{n-3}
\frac{\big|X_{\frac {k+2}{n}} -2X_{\frac {k+1}{n}}+X_{\frac {k}{n}}+X_{\frac {k+3}{n}} -2X_{\frac {k+2}{n}}+X_{\frac {k+1}{n}} \big |}
{\big |X_{\frac {k+2}{n}} -2X_{\frac {k+1}{n}}+X_{\frac {k}{n}} \big | +
\big |X_{\frac {k+3}{n}} -2X_{\frac {k+2}{n}}+X_{\frac {k+1}{n}} \big |}.-0.5174\Big ).
$$
In the R language, if {\tt X} is the vector $\big (X_{\frac 1 n},X_{\frac 2 n},\cdots, X_1 \big)$,
\begin{eqnarray*}
&&\hspace{-0.5cm} \widehat H^{(2)}_n\simeq {\tt  (mean(abs(diff(diff(X[-1]))+diff(diff(X[-length(X)]))) } \\
&& \hspace{3cm}
{\tt /(abs(diff(diff(X[-1])))+abs(diff(diff(X[-length(X)])))))-0.5174)/0.1468}.
\end{eqnarray*}
Therefore its computation is very fast and does not require any tuning parameters such as the scales for estimators based on quadratic variations or wavelet coefficients. The convergence rate of our estimator is $\sqrt n$ as for the parametric Whittle or the generalized quadratic variation estimators and hence it is  more accurate than most of other well-known semi-parametric estimators (log-periodogram, local Whittle or wavelet based estimators).

\smallskip

Estimators of the form (\ref{R_p}) can also be applied to discrete
time (sequences) instead of continuous time processes (functions).
For instance Surgailis {\it et al.}~(2008) extended the statistic
$R^{2,n}(X)$ to discrete time processes and used
it to test for $I(d)$ behavior $(-1/2< d< 5/4)$ of
observed time series. Vai\v ciulis~(2009) considered estimation  of the tail index of i.i.d. observations
using an increment ratio statistic.

\begin{remark}\label{level} The referee noted that the IR-roughness  might be connected to the level crossing index
(see Feuerverger {\it et al.}, 1994). To our surprise, such a connection indeed exists as explained below. Let
$Y_n(t), \, t \in [0,1-\frac 1 n]$ be the linear interpolation of the ``differenced'' sequence $\Delta^{1,n}_j X = X(\frac {j+1} n) - X(\frac j n), \, j=0,1, \cdots, n-1$:
$$
Y_n(t) = n \bigg[ (\frac{j+1} n - t)\Delta^{1,n}_j X  +  (t - \frac j n) \Delta^{1,n}_{j+1} X \bigg],  \qquad t \in [\frac j n, \frac{j+1} n),
$$
$j=0,1, \cdots, n-2$. Then, using Figure \ref{Figure3new} as a proof,
\begin{eqnarray} \label{IRgraph}
R^{1,n}(X)&=&\frac n {n-1} \sum_{j=0}^{n-2} \Big|{\rm meas}\Big\{ t\in [\frac j n, \frac{j+1} n): Y_n(t) > 0 \Big\}
- {\rm meas}\Big\{ t\in [\frac j n, \frac{j+1} n): Y_n(t) < 0 \Big\} \Big|~~~\\
&=&\frac n {n-1}  \sum_{j=0}^{n-2} \Big|\int_{\frac j n}^{\frac {j+1}n} (\1(Y_n(t) >0) - \1(Y_n(t)<0)) \d t \Big|. \nonumber
\end{eqnarray}
\begin{center}
\begin{figure}[h]
\begin{center}
\includegraphics[width=10 cm,height=15cm]{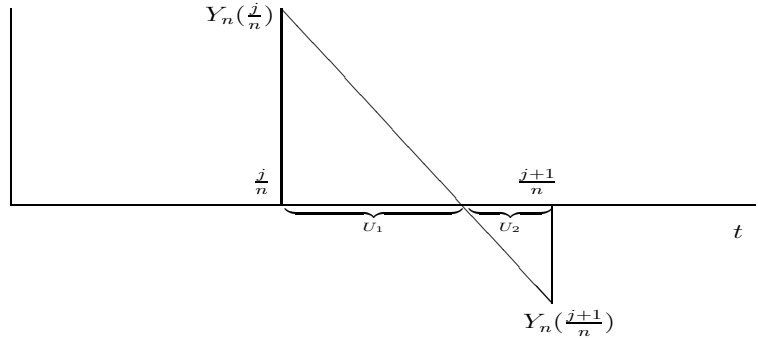}
\vspace{-10.5cm}
\caption{The proof of (\ref{IRgraph}):  follows by \
$\frac{|Y_n(\frac{j}{n}) + Y_n(\frac{j+1}{n})|}{|Y_n(\frac{j}{n})| + |Y_n(\frac{j+1}{n})|} = n |U_1 - U_2|$.  } \label{Figure3new}
\end{center}
\end{figure}
\end{center}

\noindent Let $\psi (x_1,x_2) :=  |x_1+x_2|/(|x_1|+|x_2|), \  \psi_0(x_1,x_2) := \1 (x_1 x_2 \ge 0) $.  Clearly, the two quantities
$1- \psi\big(Y_n(\frac{j}{n}),Y_n(\frac{j+1}{n})\big)$ and $1- \psi_0\big(Y_n(\frac{j}{n}),Y_n(\frac{j+1}{n})\big)$ both are strictly positive  if and only if $Y_n $ crosses the zero level in the interval $[\frac j n, \frac{j+1} n)$ but the former quantity measures not only the fact  but also the ``depth'' of the crossing so that  $1- \psi\big(Y_n(\frac{j}{n}),Y_n(\frac{j+1}{n})\big) $ attains its maximal value 1  in the case of a ``perfect'' crossing
in the middle of the interval $[\frac j n, \frac{j+1} n)$  (see Figure \ref{Figure3new}).\\
It seems that similar asymptotic results can be obtained  for  $R^{p,n}_0(X) := \frac{1}{n-p} \sum_{k=0}^{n-p-1}
\psi_0\big(\Delta^{p,n}_k X, \Delta^{p,n}_{k+1} X \big) $ measuring the number of zero crossings of the increment
sequence $\Delta^{p,n}_k X, k = 0,1, \cdots, n-p$ and other similar statistics obtained by replacing the functions $\psi $ or $\psi_0$ by other scaling invariant functions.  Let us note that $ R^{1,n}_0(X) $ is related to the zero-crossings' counting statistic studied in Ho and Sun (1987) for stationary  Gaussian time series. Also note that
the Hermite rank of $\psi_0$  is $2$ and that the corresponding limit function
$\lambda_0 (r) = \frac 1 \pi \arccos (-r) $ is strictly increasing on the interval $(-1,1)$ similarly as the function $\lambda (r) $ in (\ref{lambda_0}).
On the other hand, while the statistic $R^{p,n}_0(X)$ is certainly of interest, the statistic $R^{p,n}(X) $ seems preferable to it for the reasons
explained above.
In particularly, in the case of symmetric Lévy processes $X$ with independent increments studied in Section~\ref{Levy}, the latter statistic leads to an estimator of the fractional index while the former statistic can be easily shown to  converge to $1/2$.
\end{remark}
\smallskip

The paper is organized as follows. Section \ref{general} discusses some
general (consistency) properties of the estimators $R^{p,n}(X)$.
Section \ref{diffu} deals with the case when $X$ is a diffusion. The case of Gaussian processes $X$ is
considered in Section \ref{gauss} while the case of Lévy processes is studied in Section \ref{Levy}.
Section \ref{proofs} contains proofs and other derivations.

Below, we write $C$ for generic constants which may change from line to line.

\section{Some asymptotic results} \label{general}

The definition of $R^{p,n} f $ in (\ref{R_p}) can be extended to more general increments (the so-called generalized
variations).
Consider a filter $a:=(a_0,\cdots,a_q) \in \R^{q+1}$ such that there exists $p  \in \N,  p\le q$ satisfying
\begin{eqnarray}\label{moment}
\sum_{\ell=0}^q \ell^i a_\ell=0~~\mbox{for}~ i=0,\cdots,p-1~~~\mbox{and}~~~\sum_{\ell=0}^q \ell^{p} a_\ell \neq 0.
\end{eqnarray}
The class of such filters will be denoted  ${\cal A}(p,q)$.
For $n \in \N^* := \{1,2, \cdots \}$ and  a function $f = (f(t), t \in [0,1])$, define
the generalized variations of $f$ by
\begin{equation} \label{genvar}
\Delta^{a,n}_j f \ := \ \sum_{\ell=0}^q a_\ell  f(\frac{j+\ell}{n}), \qquad j=0,1, \cdots, n-q.
\end{equation}
A particular case of (\ref{genvar}) corresponding to $q=p\ge 1, \,  a_\ell = (-1)^{p-\ell} {p \choose \ell} $ is the $p-$order increment
$\Delta^{p,n}_j f $ in (\ref{Delta-p}). For a filter $a\in {\cal A}(p,q)$, let
\begin{eqnarray}
R^{a,n} (f) &:=&\frac{1}{n-q} \sum_{k=0}^{n-q-1}
\frac{\big|\Delta^{a,n}_k f +  \Delta^{a,n}_{k+1} f \big|}
{|\Delta^{a,n}_k f| + |\Delta^{a,n}_{k+1} f| }. \label{R_a}
\end{eqnarray}

It is easy to prove that $R^{1,n}(f) \limiten 1 $ if $f$ is continuously
differentiable on $[0,1]$ and the derivative $f'$ does not vanish
on $[0,1]$ except maybe for a finite number of points. Moreover, it is obvious that $R^{1,n}(f)  = 1$ if $f$ is monotone on
$[0,1]$: the IR-roughness of a monotone function is the same as of a smooth function,  which is
not  surprising since a similar fact holds for other measures of roughness such as the
$p-$variation index or the Hausdorff dimension.\\
We conjecture that $R^{p,n}(f) \to 1 $ and $R^{a,n}(f) \to 1$ for any $q\ge p \ge 1, \ a \in {\cal A}(p,q)$ and $f: [0,1] \to {\R} $ which is $(p-1)$
times differentiable and the derivative $f^{(p-1)}$ has bounded
variation on $[0,1]$ with the support $ {\rm supp}(f^{(p-1)}) = [0,1]$. However, we can prove a weaker result.

\begin{proposition} \label{derive}
Let $f$ be $(p-1)$-times continuously differentiable $(p\ge 1)$ with
$f^{(p-1)}$ being absolutely continuous on $[0,1]$ having the
Radom--Nikodym derivative
$g = \left(f^{(p-1)}\right)^\prime $.
Assume that
$g \ne 0 $ a.e. in $[0,1]$.
Then $R^{p,n}(f) \limiten 1$ and $R^{a,n}(f) \limiten 1 $ for any $a \in {\cal A}(p,q), \, q \ge p$.
\end{proposition}

\noindent {\it Proof.} We restrict the proof to the case  $p=2$
since the general case is analogous. Using summation by parts, we can rewrite $\Delta^{a,n}_j f $ as
\begin{eqnarray}  \label{aDelta}
\Delta^{a,n}_j f
&=&\sum_{i=0}^q b_i \, \Delta_{i+j}^{2,n} f,
\end{eqnarray}
where $b_i := \sum_{k=0}^i \sum_{\ell = 0}^k a_\ell, \, i=0,1, \cdots, q, \ b_{q-1} = b_q = 0 $ and
$ \bar b := \sum_{i=0}^q b_i = \frac 1 2 \sum_{i=1}^q i^2 a_i \ne 0 $ in view of the assumption
$a \in {\cal A}(2,q)$.

Assume $n$ is large enough
and for a given $t\in (0,1)$, let $k_n(t) \in \{0, \cdots, n-2\}$
be chosen so that $t\in [k_n(t)/n, (k_n(t)+1)/n)$, therefore
$k_n(t)=[nt]-1$. We claim that for a.e. $t \in (0,1)$
\begin{equation}
\lim_{n\to \infty} n^2 \Delta^{a,n}_{k_n(t)} f  \ = \  \bar b \, g(t),
\quad \lim_{n\to \infty} n^2 \Delta^{a,n}_{k_n(t)+1} f \ = \  \bar b \, g(t).
\label{limits}
\end{equation}
Using the fact that the function $(x_1,x_2) \mapsto \frac{|x_1 +
x_2|}{|x_1|+|x_2|} $ is continuous on ${\R}^2 \backslash \{(0,0)\}
$, we obtain
\begin{eqnarray}
h^{a,n}(t) &:=&  \frac{\left| n^2 \Delta^{a,n}_{k_n(t)} f
+ n^2 \Delta^{a,n}_{k_n(t)+1} f  \right|}
{\left| n^2 \Delta^{a,n}_{k_n(t)} f\right|
+\left|n^2 \Delta^{a,n}_{k_n(t)+1} f\right| } \
\limiten   \  \frac{|\bar b \, g(t) + \bar b \, g(t)|}{|\bar b \,g(t)| + |\bar b \, g(t)|} \  =  \ 1  \label{h_2n}
\end{eqnarray}
for a.e. $t \in (0,1)$, where we used the fact that $\bar b \, g(t) \neq 0$
a.e. Since for $n \geq q$, $R^{a,n}(f)$ can be written as $
R^{a,n}(f) = \frac n {n-q} \, \int_0^1 h^{a,n}(t) \d t$, relation
$R^{a,n}(f) \limiten 1$ follows by the dominated convergence
theorem and  the fact that $0\le h^{a,n}(t) \le 1 $.

Relations (\ref{limits}) can be proved using the Lebesgue--Vitali
theorem (see Shilov and Gurevich, 1967, Ch.~4, \S 10, Theorem 1),  as
follows.  Consider the signed measure $\mu$ on Borel subsets of
$[0,1/2]^2 $ given by
$$
\mu(A) = \int_A g(x_1+x_2) \d x_1 \d x_2.
$$
Note $\Delta^{2,n}_{k} f = \mu ((k/2n, (k+2)/2n] \times (k/2n,
(k+2)/2n]), \ (k=0, \cdots, n-2)$.  Since rectangles $[x_1, x_1 +
h]\times [x_2, x_2+h], \  (0\le x_i <x_i + h \le 1/2, i =1,2)$ form a
Vitali system on $[0,1/2]^2$, the above mentioned Lebesgue--Vitali
theorem implies that
\begin{equation} \label{limits1}
\phi_n(t_1,t_2) :=  n^2 \mu
\left(\mbox{$\left(\frac{k_{n}(t_1)}{2n},
\frac{k_{n}(t_1)+2}{2n}\right] \times  \left( \frac{k_{n}(t_2)}{2n},
\frac{k_{n}(t_2)+2}{2n}\right]$}\right) \limiten
g\left(\mbox{$\frac{t_1+t_2}{2}$}\right)
\end{equation}
a.e. in $[0,1]^2 $. Taking into account the form of the measure $\mu$ and the limiting
function in (\ref{limits1}), it follows the convergence $ n^2
\Delta^{2,n}_{k_n(t)} f = \phi_n(t,t) \limiten g(t) $ a.e. on
$[0,1]$. 
Next, for any fixed $i = 0,1, \cdots $,
the sequence of rectangles $(\frac{k_n(t_1) +i}{2n}, \frac{k_n(t_1)+i+2}{2n}] \times
(\frac{k_n(t_2) +i}{2n}, \frac{k_n(t_2)+i+2}{2n}], \ n=1,2, \cdots $ is regularly contracting to $(t_1,t_2) \in (0,1)^2 $ in the sense
of (Shilov and Gurevich, 1967, Ch.~4, \S 10). Hence, using the lemma on  p.~214 of the above monograph, it follows that
$n^2 \mu
\big( \big(\frac{k_{n}(t_1)+i}{2n},
\frac{k_{n}(t_1)+i+2}{2n}\big] \times  \big( \frac{k_{n}(t_2)+i}{2n},
\frac{k_{n}(t_2)+i+2}{2n}\big]\big) \limiten
g\big(\frac{t_1+t_2}{2}\big) $ \ a.e. in $[0,1]^2 $, implying
$$
n^2
\Delta^{2,n}_{k_n(t)+i} f  \limiten g(t)  \qquad \text{a.e. on}\
[0,1], \quad \text{for any } \ i=0,1, \cdots.
$$
Together with (\ref{aDelta}), this proves (\ref{limits}) and the proposition. \hfill $\Box$ ~\\

Let us turn now to the case when $f(t) = X_t,\, t\in [0,1]$ is a
random process. Now and in all the sequel, $R^{p,n}(X), \, R^{a,n}(X)$ are
denoted $R^{p,n}, \, R^{a,n}$, respectively.  Below we formulate a general condition for
the convergence of $R^{p,n}$ and $R^{a,n}$ to a deterministic limit.

\bigskip

\noindent {\bf Assumption (A):} \  For a.e. pairs $(t_1, t_2) \in
(0,1)^2, t_1 \neq t_2$, for $i=1,2$ there exist:
\begin{enumerate}
\item [(i)] normalizations $A^{(t_i)}(\delta) \to \infty \  \
(\delta \to 0)$,
\item [(ii)](mutually) independent random processes
$Y^{(t_i)}  = (Y^{(t_i)}(\tau), \tau \in [0,1])$,
\end{enumerate}
such that for $\delta \to 0, s_1 \to t_1, s_2 \to t_2$
\begin{equation}
~~~~~~~\left( A^{(t_1)}(\delta) (X_{s_1+\delta \tau}
- X_{s_1}), A^{(t_2)}(\delta)(X_{s_2+\delta \tau} - X_{s_2})\right)
\ \limitefdd \  (Y^{(t_1)}(\tau), Y^{(t_2)}(\tau)). \label{double}
\end{equation}

\begin{remark}
Relation (\ref{double}) implies the existence of a joint small
scale limit $(Y^{(t_1)}, Y^{(t_2)})$ at a.e. pair $(t_1, t_2)\in
(0,1)$, with independent components $Y^{(t_1)},Y^{(t_2)}$.
Note Assumption (A) and Proposition \ref{prop2.3} below are very general, in the
sense that they do not assume any particular structure or
distribution of $X$.
\end{remark}

\begin{proposition}\label{prop2.3}
Let
$a =(a_0,\cdots,a_q) \in {\cal A}(p,q), \,  1\le p \le q$ be a filter and let $X$ satisfy Assumption (A).
Assume in addition that $\P ( |\Delta^a_j Y^{(t)}| >0) =1, \, j=0,1 $ for a.e. $t \in (0,1)$, where
$\Delta^a_j Z \equiv  \Delta^{a,1}_j Z =  \sum_{\ell =0}^q a_\ell Z(j+\ell)  $.  Then
\begin{eqnarray}\label{prop_L2}
\E \bigg (R^{a,n} - \int_0^1 \E \left[ \frac{
|\Delta^a_0Y^{(t)} + \Delta^a_1 Y^{(t)}|
}{|\Delta^a_0Y^{(t)}| + |\Delta^a_1 Y^{(t)}|} \right]  \d
t\bigg)^2 \limiten 0.
\end{eqnarray}
\end{proposition}

\noindent {\it Proof.} The statement follows from
\begin{equation} \label{conv_R}
\E R^{a,n} \limiten \int_0^1 \E \left[ \frac{
|\Delta^a_0Y^{(t)} + \Delta^a_1 Y^{(t)}|
}{|\Delta^a_0Y^{(t)}| + |\Delta^a_1 Y^{(t)}|} \right]  \d
t~~~\mbox{and}~~~~ \E \left(R^{a,n} - \E R^{a,n}\right)^2 \limiten 0.
\end{equation}
Write $\E R^{a,n} = \frac n {n-q} \, \int_0^1 \E h^{a,n}_{X}(t) \d
t$, where (c.f. (\ref{h_2n}))
\begin{eqnarray*}
h^{a,n}_{X}(t) &:= &\frac{A^{(t)}(1/n)\left(\left|
\Delta^{a,n}_{k_n(t)} X + \Delta^{a,n}_{k_n(t)+1} X
\right|\right)} {A^{(t)}(1/n)\left(\left|  \Delta^{a,n}_{k_n(t)}
X\right| +\left|\Delta^{a,n}_{k_n(t)+1} X\right|
\right)} \\
& \limiteloin & \ \frac{ |\Delta^a_0Y^{(t)} + \Delta^a_1
Y^{(t)}| }{|\Delta^a_0Y^{(t)}| + |\Delta^a_1 Y^{(t)}|} \ =: \ h^a_Y (t),
\end{eqnarray*}
for a.e. $t\in (0,1)$, according to Assumption (A) and the
continuous mapping theorem. Whence and from $0\le h^{a,n}_{X}
\le 1$ using  the Lebesgue dominated convergence theorem the first relation in (\ref{conv_R}) follows. Moreover,
\begin{eqnarray*}
\E (R^{a,n} - \E R^{a,n})^2&=&\Big ( \frac n {n-q}\Big )^2 \,
\int_0^1\int_0^1 \E \big [ h^{a,n}_{X}(t) h^{a,n}_{X}(t')\big ] \d
t \d t'-\big (\E R^{a,n}\big )^2,
\end{eqnarray*}
and with the same arguments as previously and the independence
of $Y^{(t)}$ and $Y^{(t')}$ when $t\neq t'$,
$$
\E \big [ h^{a,n}_{X}(t) h^{a,n}_{X}(t') \big ] \limiten \E \bigg[
\frac{ |\Delta^a_0 Y^{(t)} + \Delta^a_1 Y^{(t)}|
}{|\Delta^a_0Y^{(t)}| + |\Delta^a_1 Y^{(t)}|} \cdot \frac{
|\Delta^a_0Y^{(t')} + \Delta^a_1 Y^{(t')}| }{|\Delta^a_0Y^{(t')}|
+ |\Delta^a_1 Y^{(t')}|}\bigg] \  =  \   \E h^{a}_{Y}(t) \cdot \E
h^{a}_{Y}(t')
$$
and therefore $\Big ( \frac n {n-q}\Big )^2 \, \int_0^1\int_0^1 \E
\big [ h^{a,n}_{X}(t)h^{a,n}_{X}(t')\big ] \d t \d t'- \Big
(\int_0^1 \E h^{a,n}_{X}(t) \d t\Big )^2 \limiten 0$,  thereby
proving  the second relation in (\ref{conv_R}) and  the
proposition. \hfill $\Box$\\
~\\
The following easy but interesting corollary can also be added to this result. It proves
that smooth multiplicative or additive trends do no change the $\L^2$-asymptotic behavior of  $R^{a,n}$. Let ${\cal C}^p[0,1] $ denote
the class of all $p-$times continuously differentiable functions on $[0,1]$.

\begin{corollary} \label{trend}
Let $a \in {\cal A}(p,q)$ and $X$ satisfy the conditions of Proposition \ref{prop2.3} with  $A^{(t)}(\delta)=O(\delta^{-1})$ $(\delta\to 0)$ for each $t \in [0,1]$. Assume that  $\alpha \in {\cal C}^1[0,1], \, \beta\in {\cal C}^p[0,1]$,  $\inf_{t\in [0,1]} \alpha(t)> 0$   and $\sup_{t \in [0,1]} |X(t)| < \infty $ a.s.
Define $Z$ such that
$Z_t=\alpha(t)\, X_t+\beta(t), \ t \in [0,1]$. Then (\ref{prop_L2})
holds with $R^{p,n}= R^{p,n}(X)$ replaced by $R^{p,n}(Z)$.
\end{corollary}

\noindent {\it Proof.} We consider here $p\geq 2$ but the case $p=1$ can be easily obtained. With
$\alpha(\frac {k+j} n )=\alpha(\frac k n)+ \frac {j} n \alpha'(\frac{k}{n}) \big (1 + o(1)\big )$ and
$\Delta^{a,n}_k \beta = O\big(\frac{1}{n^2}\big)$, for a.e. $t \in (0,1)$ and
$k  = k_n(t)$ as defined in the proof of Proposition \ref{derive}, we deduce that
\begin{eqnarray*}
\Delta^{a,n}_{k} Z&=& \alpha\big(\frac{k}{n}\big)\Delta^{a,n}_{k} X+ \frac 1 n  \alpha'\big(\frac{k}{n}\big) \Delta^{a',n}_{k}X +
o\big(\frac 1 n\big)\sup_{t \in (0,1)} |X(t)|  +  O\big(\frac{1}{n^2}\big),
\end{eqnarray*}
with $a'=(ja_j)_{0\leq j \leq q} \in  {\cal A}(p-1,q)$.
 Therefore,
\begin{eqnarray*}
A^{(t)}(\frac 1 n )\Delta^{a,n}_{k_n(t)} Z&=&\alpha\big(\frac{k_n(t)}{n}\big)A^{(t)}\big(\frac 1 n \big)\Delta^{a,n}_{k_n(t)} X + \frac 1 n O_p\big(A^{(t)}\big(\frac 1 n \big)\Delta^{a',n}_{k_n(t)} X \big)
+ o_p\big( \frac 1 n A^{(t)}(\frac 1 n )\big)\\
&\limiteloin& \alpha(t)\Delta^{a}_{0} Y^{(t)}.
\end{eqnarray*}
In a similar way, for a.e. pairs $(t, t') \in (0,1)^2, \, t\ne t'$,  we can verify the joint convergence in distribution of  r.v.'s
$A^{(t)}(\frac 1 n )\Delta^{a,n}_{k_n(t)} Z, $  $A^{(t)}(\frac 1 n )\Delta^{a,n}_{k_n(t)+1} Z, $ $ A^{(t')}(\frac 1 n )\Delta^{a,n}_{k_n(t')} Z,
A^{(t')}(\frac 1 n )\Delta^{a,n}_{k_n(t')+1} Z$ to the limiting  r.v.'s
$\alpha(t)\Delta^{a}_{0} Y^{(t)},
\alpha(t)\Delta^{a}_{1} Y^{(t)}, \alpha(t')\Delta^{a}_{0} Y^{(t')}, \alpha(t')\Delta^{a}_{1} Y^{(t')}. $  Now, the statement of the
corollary follows by the argument at the end of the proof of Proposition \ref{prop2.3}.\hfill $\Box$\\

\begin{remark}\label{final} By definition, the statistics $ R^{p,n} $ and $  R^{a,n} $ for $a \in {\cal A}(p,q), \,  1\le p \le q$
are invariant with respect to additive polynomial trends of order less than $p$; in particular, $R^{3,n}$ is insensitive to a quadratic trend while
$R^{2,n}$ does not have this property. On the other hand, Corollary \ref{trend} 
(see also Proposition \ref{trends}) says that under weak additional conditions on $X$, any sufficiently smooth additive or multiplicative trends do not affect the limit of $R^{p,n}$
as soon as $p\geq 1$.
In the important special case when the limit process $Y^{(t)} = B_H$ in Assumption (A) and (\ref{prop_L2}) is a fractional Brownian motion
with parameter $H \in (0,1)$ independent of $t$, the statistic $ R^{p,n}$ converges in mean square to the expectation
$\E \frac{|\Delta^p_0 B_H + \Delta^p_1 B_H|}{|\Delta^p_0 B_H| + |\Delta^p_1 B_H|} = \lambda (\rho_p(H)) $, c.f. (\ref{lambda_p})-(\ref{rho_p}).
Numerical computations show that the correlation coefficient $ \rho_p(H) $ is a monotone function of $H$ for any $p\ge 1 $ and tends to constant value $-1$ on the interval $(0,1)$ as $p $ increases. Therefore, for larger values of $p$, the range of $\lambda (\rho_p(H)) $ is rather small and
$ R^{p,n}$ seems less capable to estimate $H$. A final reason for our concentrating
on the ``lower-order''  statistics  $ R^{p,n}, \ p=1,2$ in the rest of the paper is the fact that $R^{2,n} $ satisfies the central limit theorem
in (\ref{CLT_fBm}) on the whole interval $H \in (0,1)$.
\end{remark}

\section{Diffusions}\label{diffu}

Let
\begin{equation}
X_t \ = \  X_0 +  \int_0^t a_s\d B(s) + \int_0^t b_s \d s, \qquad
t\in [0,1]  \label{ito}
\end{equation}
be a diffusion (or It\^o's) process on ${\R}$.  In (\ref{ito}), we
assume the existence of a right-continuous filtration ${\mathcal F}
= ({\mathcal F}_t, t \in [0,1])$, a standard Brownian motion $B$
adapted to ${\mathcal F}$; moreover,  $a_s, b_s, s \in [0,1]$ are
adapted random functions satisfying $\int_0^1 |b_s| {\d}s < \infty,
\ \int_0^1 a^2_s {\d}s < \infty $ a.s., and $X_0$ is a ${\mathcal
F}_0-$measurable r.v. Write $\E_t [\cdot ] = \E [\cdot |{\cal F}_t]
$ for the conditional expectation. Let $\Lambda_1(1/2) =
\lambda(\rho_1 (1/2)) \simeq 0.7206$ and $\Lambda_2(1/2) =
\lambda(\rho_2 (1/2)) \simeq 0.5881$. The proof of the following Lemma \ref{Lem32}
is given in Annexe.

\begin{lemma} \label{Lem32}
Let $\psi(x_1,x_2) :=
|x_1+x_2|/(|x_1|+|x_2|) \ (x_1,x_2 \in {\R})$, and let
 $Z_i, i=1,2 $ be independent
${\mathcal N}(0,1)$ r.v.'s. Then for any random variables $\xi_1,
\xi_2$,
\begin{equation}
\left| \E \psi(Z_1+\xi_1,Z_2+\xi_2) - \E \psi (Z_1,Z_2)\right| \ \le
\  20 \max_{i=1,2} \left(\E \xi_i^2 \right)^{1/3}. \label{Z_bdd}
\end{equation}
\end{lemma}

\smallskip

\begin{theorem} \label{Theo31} Assume the following
conditions: there exist random variables $K_1, K_2 $ such that $0<
K_i < \infty $ a.s., and such that, for any sufficiently small $h>0$
and any $0\le t < t+ h \le 1$, the following inequalities hold,
a.s.:
\begin{eqnarray}
|a_t| &\ge & K_1, \qquad \E_t b^2_{t+h} \ \le  \  K_2   \quad
\text{and} \quad  \E_t (a_{t+h}- a_t)^2 \  \le  \ K_2h.  \label{K_3}
\end{eqnarray}
Then
\begin{equation}
R^{p,n} - \Lambda_p(1/2)  \ = \  O(n^{-1/3})  \quad  \text{ a.s.}
\quad (p=1,2).
\end{equation}
\end{theorem}

\noindent {\it Proof.} We restrict the proof to the case $p=1$ since
the case $p=2$ is analogous. For notational simplicity, assume that $n $ is odd. Define
\begin{equation}
\eta_n(k) := \frac{|\Delta^{1,n}_k X + \Delta^{1,n}_{k+1}
X|}{|\Delta^{1,n}_k X| + |\Delta^{1,n}_{k+1} X|}, \quad \eta'_n(k)
:= \E_{k/n} [\eta_{n}(k)], \quad \eta''_n(k) := \eta_{n}(k)-
\eta'_n(k)
\end{equation}
and correspondingly write $R^{1,n} = R'_n + R''_n, \  R'_n :=
(n-1)^{-1}\sum_{k=0}^{n-2} \eta'_n(k), \ R''_{n1} :=
(n-1)^{-1}\sum_{k=0}^{(n-2)/2} \eta''_n(2k), $   $  R''_{n2} :=
(n-1)^{-1}\sum_{k=0}^{(n-4)/2} \eta''_n(2k+1). $ As $(\eta''_n (2k),
{\mathcal F}_{(2k+2)/n}, k=0, \dots, (n-2)/2) $ is a martingale
difference sequence, so by Burkholder's inequality,
$$
\E (R''_{n1})^8 \le C n^{-8} \left(\sum_{k=0}^{(n-2)/2} \E^{1/4}
(\eta'_n(2k))^8 \right)^4 \le Cn^{-4}
$$
and therefore
$$
\sum_{n=1}^\infty \P (|R''_{n1}|> n^{-1/3} ) \le C\sum_{n=1}^\infty
n^{8/3} n^{-4} < \infty,
$$
implying $R''_{n1} = O(n^{-1/3})$ a.s. A similar fact holds for $R''_{n2}$.
Thus, it remains to prove
\begin{equation}
R'_n - \Lambda_1(1/2)   = O(n^{-1/3})  \qquad a.s. \label{R'_n}
\end{equation}
Observe
$$
\eta'_n(k) - \Lambda_1(1/2) \ = \  \E_{k/n}
\left[\frac{|Z_1(k)+\xi_1(k) + Z_2(k)+\xi_2(k)|} {|Z_1(k) +
\xi_1(k)| + |Z_2(k)+\xi_2(k)|}\right] - \E \left[
\frac{|Z_1(k)+Z_2(k)|}{|Z_1(k)|+|Z_2(k)|}\right],
$$
where
\begin{eqnarray*}
Z_1(k)&:=&n^{1/2}\Delta^{1,n}_k B, \qquad Z_2(k) \ := \ n^{1/2} \Delta^{1,n}_{k+1} B, \\
\xi_1(k)&:=&n^{1/2}\int_{k/n}^{(k+1)/n} \Big(\frac{a_s}{a_{k/n}} -1\Big) \d B(s) + n^{1/2}\int_{k/n}^{(k+1)/n} \frac{b_s}{a_{k/n}}\d s,\\
\xi_2(k)&:=&n^{1/2}\int_{(k+1)/n}^{(k+2)/n} \Big(\frac{a_s}{a_{k/n}}
-1\Big) \d B(s) + n^{1/2}\int_{(k+1)/n}^{(k+2)/n}
\frac{b_s}{a_{k/n}}\d s.
\end{eqnarray*}
According to Lemma \ref{Lem32} above, $|\eta'_n(k) - \Lambda_1(1/2)| \le 36
\max_{i=1,2} \left(\E_{k/n} \xi^2_i(k)\right)^{1/3} $ and therefore
$$
|R'_n - \Lambda_1(1/2)| \ \le \  36 \, \max \{ (\E_{k/n}
\xi^2_i(k))^{1/3}: i=1,2, k=0,1, \dots, n-1\}.
$$
Whence, (\ref{R'_n}) follows from the following fact: there exists a
r.v. $K< \infty $, independent of $n$ and such that for any $n\ge 1,~
k=0, \dots, n-1,~ i=1,2$
\begin{equation}\label{diff}
\E_{k/n} \xi^2_i(k) \  \le \ K\, n^{-1},  \quad \text{a.s.}
\end{equation}
Indeed, using (\ref{K_3}),
\begin{eqnarray*}
\E_{k/n}\xi_1^2(k)
&=&n\int_{k/n}^{(k+1)/n}\E_{k/n}\Big(\frac{a_s}{a_{k/n}} -1\Big)^2
\d s + n \E_{k/n} \Big(\int_{k/n}^{(k+1)/n} \frac{b_s}{a_{k/n}}\d s \Big)^2 \\
&\le& n K_1^{-2} \int_{k/n}^{(k+1)/n}\E_{k/n}\Big(a_s -
a_{k/n}\Big)^2 \d s +
K_1^{-2} \int_{k/n}^{(k+1)/n} \E_{k/n} b^2_s\d s \\
&\le& K_2K_1^{-2} n^{-1}, \qquad {\rm a.s.},
\end{eqnarray*}
and the bound (\ref{diff}) for $i=2$ follows similarly. This proves
(\ref{diff}) and Theorem \ref{Theo31}, too. \hfill $\Box$

\bigskip
\noindent Let us present some examples of It\^o's processes $X$ satisfying
conditions  (\ref{K_3}).\\
~\\
\noindent {\sc Example 3.1} \  Let $(X_t,\, t\in [0,1])$ be a Markov
process satisfying a stochastic equation
\begin{equation}
X_t = x_0 + \int_0^t  a(X_s) \d B(s) + \int_0^t b(X_s) \d s,
\label{Ito_eq}
\end{equation}
where $x_0 \in {\R}$ is nonrandom,  $a(x), b(x), x\in {\R}$ are real
measurable functions and $B$ is a standard Brownian motion. Let
${\mathcal F}_t := \sigma \{B(s), s\le t\}, 0\le t \le 1$ be the
natural  filtration. Assume that
\begin{equation}
|a(x)- a(y)| \le K|x-y|, \qquad  |b(x)-b(y)| \le K|x-y| \qquad (x, y
\in {\R})  \label{Lip}
\end{equation}
for some constant $K<\infty$.  Then equation (\ref{Ito_eq}) admits a
unique adapted solution; see e.g. Gikhman and Skorohod (1969). Let
$a_t = a(X_t), b_t = b(X_t)$. Assume in addition that $|a(x)| \ge
K_1, \ (x\in {\R})$ for some nonrandom constant $K_1>0$. Then the
first inequality in (\ref{K_3}) is trivially satisfied; moreover,
the second and third relations in (\ref{K_3}) are also
satisfied, with $K_2 = C(1+ \sup_{0\le t\le 1} X^2_t) < \infty$ and $K_3
= C$, where $C$ is nonrandom and depends on the constant $K$ in
(\ref{Lip}) only.

\bigskip

\noindent {\sc Example 3.2} \  Let $X_t := g(t, B(t)), $ where $B$
is a standard Brownian motion and $g(t,x)$ is a (jointly) continuous
function on $[0,1]\times {\R}$, having continuous partial
derivatives $g_t(t,x) := \partial g(t,x)/\partial t, $ $g_x(t,x) :=
\partial g(t,x)/\partial x, \ g_{xx}(t,x) =
\partial^2 g(t,x)/\partial x^2 $. By It\^o's lemma,
$$
\d X_t = g_x(t,B(t)) \d B(t) + \left( g_t(t,B(t)) + \frac{1}{2}
g_{xx}(t,B(t)) \right) \d t,
$$
so that $X$ admits the representation (\ref{ito}) with $a_t = g_x(t,
B(t)), b_t =  g_t(t,B(t)) + \frac{1}{2} g_{xx}(t,B(t)) $ and the
same  filtration as in the previous example. Assume that
$$
|g_x(t,x)|\ge K_1, \qquad |g_x(s,y) - g_x(t,x)| \le K(|s-t|^{1/2} +
|y-x|),
$$
for all $(t,x), (s,y) \in [0,1]\times {\R})$ and some constants $0<
K_1, K < \infty $. Then $X$ satisfies the conditions in (\ref{K_3}).

\section{Gaussian processes}\label{gauss}

\subsection{Assumptions}

Let $X = (X_t,\, t\in [0,1])$ be a Gaussian process, with
zero mean. Without loss of generality, assume $X_0 =0$. Define
$\sigma^2_{p,n}(k)$, the variance of $\Delta^{p,n}_k X$, and
$\rho_{p,n}(k)$, the correlation coefficient between
$\Delta^{p,n}_{k}X $ and $\Delta^{p,n}_{k+1}X$, {\it i.e.}
\begin{eqnarray}
\sigma^2_{p,n}(k)&:=&\E \Big[ \left(\Delta^{p,n}_k X\right)^2\Big], \qquad
\rho_{p,n}(k)\ := \ \frac{\E  \big [ \Delta^{p,n}_{k} X  \   \Delta^{p,n}_{k+1}
X\big ]}{\sigma_{p,n}(k)\sigma_{p,n}(k+1)}.
\label{rho_nk1}
\end{eqnarray}
Let $B_H = (B_H(t), t\in {\R}) $ be a fractional Brownian motion
(fBm) with parameter $0<H<1$,  i.e., a Gaussian process with zero
mean and covariance such that $\E B_H(s) B_H(t)=
\frac{1}{2}\left(|t|^{2H} + |s|^{2H} - |t-s|^{2H}\right)$. Its
$p-$th order increments $\left(\Delta^p_j B_H, j \in {\Z}\right) $
form a stationary Gaussian process, for any $p\ge 1$. In
particular, the covariance function of $\Delta_j B_H \equiv
\Delta^1_j B_H = B_H(j+1)-B_H(j)$ and $\Delta^2_j B_H =
B_H(j+2)-2B_H(j+1)+ B_H(j)$ can be explicitly calculated:
\begin{eqnarray}
\E \big [ \Delta_0B_H \Delta_j B_H\big ]
&=&
2^{-1} \left(|j+1|^{2H}+|j-1|^{2H}-2|j|^{2H}\right),  \label{cov1}\\
\E \big [ \Delta^2_0B_H \Delta^2_j B_H\big ]
&=&
2^{-1} \left(-|j+2|^{2H}+4|j+1|^{2H} - 6|j|^{2H}
+4|j-1|^{2H}-|j-2|^{2H}\right). \label{cov2}
\end{eqnarray}
From Taylor expansion,
\begin{eqnarray*}
\E \big [ \Delta_0B_H \Delta_j B_H\big ]
&\sim&2H(2H-1)j^{2H-2}, \\
\E \big [ \Delta^2_0B_H \Delta^2_j B_H\big ]
&\sim&2H(2H-1)(2H-2)(2H-3)j^{2H-4},
\end{eqnarray*}
as $j \to \infty$, and therefore the first increment, $(\Delta_j B_H) $, has a summable
covariance if and only if $0< H< 3/4 $, while the second increment,
$(\Delta^2_j B_H)$, has a summable covariance for any $0<H<1$.

\bigskip

Introduce the following conditions:
\begin{enumerate}
\item [{\bf (A.1)}] There exist continuous functions $H(t)\in (0,1)$ and
$c(t)>0$ for $t \in [0,1]$ such that $\forall  j\in {\N}^*$
\begin{eqnarray}
\lim_{n\to \infty} \sup_{t \in (0,1)}  \left| \frac{  \E  \big(
X_{([nt]+j)/n} - X_{[nt]/n}  \big)^2 }{  \left(j/n\right)^{2H(t)}
} - c(t) \right| \  = \  0,~~~\mbox{with}  \label{S_j1}
\end{eqnarray}
\begin{equation}
\lim_{n\to \infty} \sup_{t \in (0,1)} \big |H(t) -  H(t + \frac 1 n) \big| \log n = 0. \label{H_cont}
\end{equation}

\item [{\bf (A.1)}$^\prime$] There exist continuous functions $H(t)\in (0,1)$ and
$c(t)>0$ for $t \in [0,1]$ such that $\forall  j\in {\N}^*$
\begin{eqnarray}
\lim_{n\to \infty} \sup_{t \in (0,1)} \sqrt{n} \left| \frac{  \E
\big( X_{([nt]+j)/n} - X_{[nt]/n}  \big)^2 }{
\left(j/n\right)^{2H(t)}  } - c(t) \right| \  = \  0,~~~\mbox{with}
\label{S_j11}
\end{eqnarray}
\begin{equation}
~~~~~~~\lim_{n\to \infty} \sup_{t \in (0,1)} \big |H(t) -  H(t + \frac 1 n)\big | \sqrt{n} \log n = 0
~~~\mbox{and}~~~\lim_{n\to \infty} \sup_{t \in (0,1)}
\big |c(t) -  c(t + \frac 1 n)\big |\sqrt n  = 0. \label{H_cont1}
\end{equation}

\item [{\bf (A.2)}$_p$] There exist $d>0$, $\gamma> 1/2$ and $0 \le \theta < \gamma/2$
such that for any $1\le k <j \le n$ with $n\in {\N}^*$
\begin{equation} \label{rho_dom2}
\Big |\E \big [\Delta^{p,n}_{k} X \  \Delta^{p,n}_{j} X \big ]
\Big | \ \le \ d \ \sigma_{p,n}(k) \sigma_{p,n}(j) \cdot n^{\theta} \cdot
|j-k|^{-\gamma}.
\end{equation}
\end{enumerate}

\smallskip
A straightforward application of Assumption {\bf (A.1)} (or {\bf
(A.1)}$^\prime$) implies that $\sqrt{c(t)}B_{H(t)}$ is the tangent
process of $X$ for all $t \in (0,1)$ and more precisely:
\begin{property} \label{prop31}
Assumptions {\bf (A.1)}, {\bf (A.1)}$^\prime$ respectively imply
that, for any $j\in {\Z}$ and $p=1,2$,
\begin{eqnarray}
\lim_{n\to \infty} \sup_{t \in (0,1)}  \left| \frac{ \E \big [
\Delta^{p,n}_{[nt]} X \  \Delta^{p,n}_{j + [nt]} X\big ] }
{\E  \big [ \Delta^{p,n}_0 B_{H(t)} \  \Delta^{p,n}_{j} B_{H(t)} \big] }
- c(t) \right| \  = \  0,   \label{R_j1} \\
\lim_{n\to \infty} \sqrt{n} \sup_{t \in (0,1)}  \left| \frac{ \E \big [
\Delta^{p,n}_{[nt]} X \  \Delta^{p,n}_{j + [nt]} X\big ] }
{\E  \big [ \Delta^{p,n}_0 B_{H(t)} \  \Delta^{p,n}_{j} B_{H(t)} \big] }
- c(t) \right| \  = \  0.  \label{R_j11}
\end{eqnarray}
Moreover, for any $t\in (0,1)$ and $p=1,2$
\begin{eqnarray*}
\left(n^{H(t)}\Delta^{p,n}_{j + [nt]} X\right)_{j\in \Z} \limitefddn
\left(\sqrt{c(t)}\,\Delta^p_j B_{H(t)}\right)_{j \in
\Z}.\label{fdd_Delta}
\end{eqnarray*}

\smallskip

\end{property}

Assumption {\bf (A.1)} can be characterized as {\it uniform local
self-similarity of $(X_t)$} (the uniformity refers to the
supremum over $t\in (0,1)$ in (\ref{S_j1})).  Note that for $X$
having stationary increments and variogram $V(t) = \E X_t^2$,
Assumption {\bf (A.1)}  reduces to $V(t) \sim c t^{2H} \  (c>0,
0<H < 1)$. For $j=0, 1$, relation (\ref{R_j1}) implies that for any $t \in (0,1)$, the
variance and the $(1/n)-$lag correlation coefficient of
$\Delta^{p,n}_{[nt]} X$ satisfy the following relations:
\begin{eqnarray}
~~~~~~~\sigma^2_{1,n}([nt])&\simn&c(t)\sigma_1^2(H(t)) =  c(t) \,  \E \big
[ (\Delta_0 B_{H(t)})^2 \big ] =
c(t)\Big ( \frac 1 n\Big )^{2H(t)},  \label{sigma1} \\
~~~~~~~\rho_{1,n}([nt])&\limiten&\rho_1 (H(t)) = \ {\rm corr}(B_{H(t)}(1),
B_{H(t)}(2) - B_{H(t)}(1)) = 2^{2H(t)-1} - 1, \label{corr1} \\
~~~~~~~\sigma^2_{2,n}([nt])&\simn&c(t)\sigma^2_2(H(t)) = c(t) \,  \E \big [
(\Delta^2_0 B_{H(t)})^2 \big ] =
c(t) (4- 4^{H(t)}) \Big (\frac 1 n\Big )^{2H(t)},  \label{sigma2} \\
\rho_{2,n}([nt])&\limiten&\rho_2 (H(t)) = \ {\rm corr}(B_{H(t)}(2) -
2B_{H(t)},
B_{H(t)}(3) - 2B_{H(t)}(2) + B_{H(t)}(1))\nonumber \\
&=&\frac{-3^{2H(t)}+2^{2H(t)+2} -7}{8-2^{2H(t)+1}} \label{corr2}
\end{eqnarray}
(see (\ref{rho_1})). Moreover, relations
(\ref{sigma1})--(\ref{corr2}) hold uniformly in $t\in (0,1)$.
Condition (\ref{H_cont}) is a technical condition which implies (and
is "almost equivalent''  to) the continuity of the function $t \to
H(t)$. Assumption ${\bf (A.1)}'$ is a sharper convergence condition
than Assumption {\bf (A.1)} required for establishing central limit
theorems.

\smallskip

Condition (\ref{rho_dom2}) specifies a nonasymptotic
inequality satisfied by the correlation of increments
$\Delta^{p,n}_{k}X$. The particular case of stationary processes
allows to better understand this point. Indeed, if $(X_t)$ has
stationary increments, the covariance of the stationary process
$(\Delta^{p,n}_{k}X, k \in {\Z})$ is completely determined by the
variogram $V(t)$, e.g.
\begin{eqnarray}
\E \big [\Delta^{1,n}_{k} X \  \Delta^{1,n}_{j} X \big ]= \frac 1
2 \Big\{ V \Big ( \frac {k-j+1} n  \Big )  + V \Big (\frac { k-j-1} n
\Big  ) - 2V \Big (\frac {k-j}n \Big )\Big \}. \label{vario}
\end{eqnarray}
In the ``most regular'' case, when $X=B_H$ is a fBm and therefore
$V(t)=t^{2H}$, it is easy to check that assumption {\bf (A.2)}$_2$
holds with $\theta = 0$ and $\gamma = 4-2H>2 \  (0<H<1)$, while
{\bf (A.2)}$_1$ with $\theta =0,  \gamma = 2-2H$ is equivalent to
$H < 3/4$ because of the requirement $\gamma > 1/2$.  However, for
$X=B_H$, {\bf (A.2)}$_1$ holds with appropriate $\theta>0$ in the
wider region $0< H< 7/8$, by choosing $\theta < 2-2H$ arbitrarily
close to $2-2H$ and then $\gamma <2-2H + \theta$ arbitrary  close
to $4-4H$. A similar choice of parameters $\theta $ and $\gamma $
allows to satisfy {\bf (A.2)}$_p$ for more general $X$ with
stationary increments and variogram $V(t) \sim c t^{2H} \ (t \to
0)$, under additional  regularity conditions on $V(t)$ (see
below).

\begin{property} \label{prop_V}
Let $X$ have stationary increments and variogram $V(t) \sim
ct^{2H} \  (t \to 0)$, with $c>0,~ H\in (0,1)$.

\noindent (i) Assume, in addition, that  $0<H<7/8$ and
$|V''(t)| \le Ct^{-\kappa} \ (0<t<1)$, for some $C>0 $ and
$4-4H> \kappa \ge 2 - 2H, \kappa > 1/2$. Then assumption {\bf (A.2)}$_1$ holds.

\smallskip

\noindent (ii) Assume, in addition, that  $|V^{(4)}(t)| \le Ct^{-\kappa}\  (0<t<1)$, for some $C>0 $ and
$8-4H> \kappa \ge 4 - 2H$. Then assumption {\bf (A.2)}$_2$ holds.

\end{property}

The following property provides a sufficient condition for {\bf
(A.2)}$_p$ in spectral terms, which does not require
differentiability of the variogram.

\begin{property} \label{prop_f}
Let $X$ be a Gaussian process having stationary increments and the
spectral representation  (see for instance Cram\`er and
Leadbetter, 1967)
\begin{equation} \label{spec_V}
X_t\ = \ \int_{\R} \big(\e^{\i t\xi} - 1\big)\, f^{1/2}(\xi) \,
W(\d\xi),~~~~\mbox{for all}~~t \in \R,
\end{equation}
where $W({\d}x) = \overline{W (-\d x)}$ is a complex-valued
Gaussian white noise with zero mean and variance $\E
\left|W({\d}x)\right|^2 = \d x$ and $f$ is a non-negative
even function called the spectral density of $X$ such that
\begin{equation} \label{cond:f}
\int_{\R} \big(1\wedge |\xi|^2\big)\, f (\xi) \, \d \xi \ < \  \infty.
\end{equation}
Moreover, assume that $f$ is differentiable on $(K, \infty)$ and
\begin{equation} \label{spec_f}
f(\xi) \sim c \, \xi^{-2H-1} \quad (\xi\to \infty),  \qquad |f'(\xi)|
\le C \, \xi^{-2H-2} \quad (\xi>K)
\end{equation}
for some constants $c, C, K >0$. Then $X$  satisfies assumption
{\bf (A.2)}$_1$ for $0<H<3/4$ and assumption {\bf (A.2)}$_2$ for \
$0<H<1$.
\end{property}

\subsection{Limit theorems}

Before establishing limit theorems for the statistics $R^{p,n}$
for Gaussian processes, recall that $\lambda$ is given in
(\ref{lambda_0}) and with $\rho_p(H)$ in (\ref{rho_1}) one has
\begin{eqnarray*}
\int_0^1 \Lambda_p(H(t)) \d t \  = \int_0^1
\lambda(\rho_p(H(t))) \d t \  =  \  \int_0^1 \E \left [
\frac{\left|\Delta^p_0 B_{H(t)} + \Delta^p_1 B_{H(t)}\right|}
{\left|\Delta^p_0 B_{H(t)}\right| + \left|\Delta^p_1
B_{H(t)}\right|}\right ] \d t.
\end{eqnarray*}
Straightforward computations show that Assumptions {\bf (A.1)} and {\bf (A.2)$_p$ } 
imply Assumption
{\bf (A)} with $A^{(t)}(\delta)=\delta^{-H(t)}, \, Y^{(t)}= \sqrt{c(t)}\, B_{H(t)}$. Therefore Proposition
\ref{prop2.3} ensures the convergence (in $\L^2$) of the statistics
$R^{p,n}$ to $\int_0^1 \Lambda_p(H(t)) \d t$. Bardet and Surgailis~(2009) proved a.s. convergence  in Theorem \ref{psgauss}, below, using a general moment
bound for functions of multivariate Gaussian processes (see Lemma \ref{lemgauss} in Section \ref{proofs}). A sketch of this proof can be found in Section \ref{proofs}.

\begin{theorem}\label{psgauss}
Let $X$ be a Gaussian process satisfying Assumptions {\bf (A.1)}
and {\bf (A.2)$_p$}. Then,
\begin{eqnarray}
R^{p,n}  \limitepsn \int_0^1 \Lambda_p(H(t)) \d t  \qquad (p=1,2). \label{R_n1}
\end{eqnarray}
\end{theorem}

\begin{corollary}\label{cor31}
Assume that $X$ is a Gaussian process having stationary
increments, whose variogram satisfies the conditions of Properties
\ref{prop_V} or  \ref{prop_f}. Then
\begin{eqnarray}\label{psgaussstat}
R^{p,n}  \limitepsn  \Lambda_p(H) \qquad (p=1,2).
\end{eqnarray}
\end{corollary}

The following Theorem \ref{TLCgauss} is also established in Bardet and Surgailis~(2009). Its proof (see a sketch of this proof in Section \ref{proofs}) uses a general central limit theorem for Gaussian subordinated nonstationary triangular
arrays (see Theorem \ref{tlcgauss} in Section \ref{proofs}). Note that the Hermite rank of $\psi(x_1,x_2)=|x_1+x_2|/(|x_1|+|x_2|)$ is $2$ and this explains the difference between the cases $p=1$ and $p=2$ in Theorem \ref{TLCgauss}: in the first case,
the inequalities in (\ref{rho_dom}) for $({\mbf Y}_n(k)) $ as defined in (\ref{Y_11})-(\ref{Y_21}) hold only if $\sup_{t \in[0,1]}H(t) <3/4$,  while
in the latter case these inequalities hold for  $0<\sup_{t \in[0,1]}H(t)<1$. A similar fact is true also   for the estimators based on generalized quadratic variations, see Istas and Lang~(1997), Coeurjolly~(2001).

\begin{theorem}\label{TLCgauss}
Let $X$ be a Gaussian process satisfying assumptions {\bf
(A.1)}$^\prime$ and {\bf (A.2)}$_p$, with  $\theta = 0$.  Moreover,
assume additionally $\displaystyle{\sup_{t \in[0,1]}H(t) <3/4}$ if
$p=1$. Then, for $p=1, 2$,
\begin{equation}
\sqrt{n}\Big (R^{p,n}  - \int_0^1 \Lambda_p(H(t)) \d t \Big )
\limiteloin {\mathcal N} \Big (0,\int_0^1 \Sigma_p(H(\tau))\d \tau
 \Big ),
\label{TLC}
\end{equation}
with $\Lambda_p(H)$ and $\Sigma_p(H)$
given in  (\ref{lambda_p}) and (\ref{Sigma_p}), respectively.
 \end{theorem}

The following proposition shows that the previous theorems are satisfied when smooth multiplicative and additive trends are considered.

\begin{proposition}\label{trends} Let $Z = (Z_t= \alpha(t)X_t+  \beta(t),\, t \in [0,1])$, where
$X = (X_t, t \in [0,1]) $ is a zero mean Gaussian process and $\alpha, \beta $ are deterministic continuous
functions on $[0,1]$ with $\inf_{t\in [0,1]} \alpha(t)>0$.

\smallskip

\noindent (i) \, Let $X$ satisfy assumptions of Theorem \ref{TLCgauss} and $\alpha \in {\cal C}^p [0,1],  \beta \in {\cal C}^p [0,1].$
Then the statement of Theorem \ref{TLCgauss} 
holds with $X$ replaced by $Z$.

\smallskip

\noindent (ii) \,  Let $X$ satisfy assumptions of Theorem \ref{psgauss} and $\alpha \in {\cal C}^1 [0,1],  \beta \in {\cal C}^1 [0,1].$
Then the statement of Theorem \ref{psgauss}
holds with $X$ replaced by $Z$.
\end{proposition}

\begin{remark}  \label{weakTLC}
A version of the central limit theorem in (\ref{TLC}) is
established in Bardet and Surgailis~(2009) with $\int_0^1 \Lambda_p(H(t)) \d t$  replaced by $\E
R^{p,n}$ under weaker assumption than {\bf (A.1)}$^\prime$ or even
{\bf (A.1)}: only properties (\ref{sigma1})-(\ref{corr1}) (for $p=1$) and (\ref{sigma2})-(\ref{corr2}) (for $p=2$),
in addition to {\bf (A.2)}$_p$ with  $\theta = 0$ ,
are required.
\end{remark}
The particular case of Gaussian processes having stationary
increments can also be studied:

\begin{corollary}\label{corr0}
Assume that $X$ is a Gaussian process having stationary increments
and there exist $c>0$, $C>0$ and $0<H<1 $ such that at least one of
the two following conditions (a), (b) hold:
\begin{eqnarray*} \label{condprop3}
&&\text{(a)} \  \text{variogram} \quad  V(t) = c\,  t^{2H}\big
(1+o(t^{1/2})\big )~~\mbox{for}~~ t \to 0~~\mbox{and}~~|V^{(2p)}(t)|
\le C
t^{2H-2p}~~\mbox{for all}~~t\in (0,1];\\
&&\text{(b)}\ \mbox{spectral density $f$ satisfies (\ref{cond:f}),
(\ref{spec_f}) and $f(\xi) = c \xi^{-2H-1}\big (1+o(\xi^{-1/2})\big
) \quad (\xi\to \infty)$.}
\end{eqnarray*}
Then:
\begin{equation}\label{TLCstat}
\sqrt{n} \big(R^{p,n} - \Lambda_p(H) \big) \limiteloin {\cal
N}\Big(0,\Sigma_p(H) \Big) \qquad \text{if}\ ~~~ \cases{p=1, \   \
0< H < 3/4, \cr p=2, \   \    0<H<1. \cr}
\end{equation}
Moreover, with the expression of $s_2^2(H)$ given in
Section~\ref{proofs},
\begin{equation}
\sqrt{n} \Big (\Lambda^{-1}_2(R^{2,n})-H \Big) \limiteloin {\cal
N}\big(0,s_2^2(H)\big ). \label{TLCstat2}
\end{equation}
\end{corollary}
Therefore, $\widehat H_n=\Lambda^{-1}_2(R^{2,n})$ is an estimator of
the parameter $H$ following a central limit theorem with a
convergence rate $\sqrt n$ under semi-parametric assumptions.
Similar results were obtained by Guyon and Leon~(1989) and Istas and
Lang~(1997) for generalized quadratic variations under less general assumptions.

\begin{remark}
In the context of Corollary \ref{corr0} and $H \in (3/4,1)$, we
expect  that $R^{1,n}$ follows a nongaussian limit distribution
with convergence rate $n^{2-2H}$.
\end{remark}

\begin{remark} \label{Hlinear}
Figure 1  exhibits that $H \mapsto \Lambda_2(H)$ is nearly
linear and is well-approximated by $.1468 H+.5174$. Consequently, $\int_0^1 \Lambda_2(H(t)) \d t \approx
.1468 \bar H+.5174$, where $\bar H=\int_0^1 H(t) \d t$ is the mean
value of the function $H(\cdot )$.
\end{remark}

Another interesting particular case of Theorem \ref{TLCgauss}
leads to a punctual estimator of the function $H(t)$ from a
localization of the statistic $R^{2,n}$. For $t_0 \in
(0,1)$ and $\alpha \in (0,1)$, define
\begin{eqnarray*}
R_\alpha^{2,n}(t_0) := \frac 1 {2 n^\alpha}
\sum_{k=[nt_0-n^\alpha]}^{[nt_0+n^\alpha]}
\frac{\left|\Delta^{2,n}_{k} X  + \Delta^{2,n}_{k+1} X\right|}
{|\Delta^{2,n}_{k}X| + |\Delta^{2,n}_{k+1}X| }. \label{R_t0}
\end{eqnarray*}
This estimator is studied in Bardet and Surgailis~(2010) and compared to the estimator based on generalized quadratic variations
discussed in Benassi {\it et al.}~(1998) and Coeurjolly~(2005).

\subsection{Examples}
Below we provide some concrete examples of Gaussian processes which
admit a fBm as the tangent process. For some examples the hypotheses
of Theorems \ref{psgauss}-\ref{TLCgauss} and the subsequent corollaries are satisfied. For
other examples, the verification of our hypotheses (in particular,
of the crucial covariance bound ${\bf (A.2)}_p $) remains an open
problem and will be discussed elsewhere.

\begin{example} {\em Fractional Brownian motion (fBm)}. As noted above,
a fBm $X= B_H$ satisfies  {\bf (A.1)}$^\prime$ as well as {\bf
(A.2)}$_1$ (for $0<H<3/4$ if $\theta=0$ and $0<H<7/8$ if $0<\theta<2-2H$ with $\theta$ arbitrary close to $2-2H$ and therefore $\gamma<2-2H+\theta$ arbitrary close to $4-4H$ may satisfy $\gamma>1/2$) and {\bf (A.2)}$_2$ (for $0<H<1$), with
$H(t) \equiv H$, $c(t)\equiv c$. Therefore, for fBm both Theorems
\ref{psgauss} (the almost sure convergence, satisfied for $0<H<7/8$ when $p=1$ and for $0<H<1$ when $p=2$) and \ref{TLCgauss} (the central limit
theorem, satisfied for $0<H<3/4$ when $p=1$ and for $0<H<1$ when $p=2$)) apply.
Obviously, a fBm also satisfies the conditions of  Corollary
\ref{corr0}. Thus, the rate of convergence of the estimator
$\Lambda_2^{-1}(R^{2,n}) =: \widehat H_n$ of $H$ is $\sqrt{n}$.
But in such a case the self-similarity property of fBm allows to use
in this case asymptotically efficient Whittle or maximum likelihood
estimators (see Fox and Taqqu, 1987, or Dahlhaus, 1989). However,
for a fBm with a continuously differentiable multiplicative and
additive trends, which leads to a semi-parametric context, the
convergence rate of $\widehat H_n$ is still $\sqrt{n}$ while parametric estimators cannot be applied.
\end{example}

\begin{example} {\em Multiscale fractional
Brownian motion} (see Bardet and
Bertrand, 2007) defined as follows: for $\ell \in \N^*$, a
$(M_\ell)$-multiscale fractional Brownian motion $X=(X_t, \, t\in
\R)$ ($(M_\ell)$- fBm for short) is a Gaussian process having
stationary increments and a spectral density $f$ such that
\begin{equation}
f(\xi) =   \frac {\sigma_j^2}{|\xi|^{2H_j+1}} \,
\1(\omega_j\le  |\xi| < \omega_{j+1})~~~\mbox{for all}~~\xi\in \R
\label{MK_fBm}
\end{equation}
with $\omega_0 := 0<\omega_1<\cdots<\omega_\ell<\omega_{\ell+1}
:=\infty, \ \sigma_i >0$ and $H_i \in \R$ for $i \in
\{0,\cdots,\ell\}$ with $H_0<1$ and $H_\ell>0$. Therefore condition
(\ref{spec_f}) of Property \ref{prop_f} is satisfied, with $K =
\omega_\ell$ and $H=H_\ell$. Moreover, the condition $f(\xi) = c
\xi^{-2H-1}\big (1+o(\xi^{-1/2})\big ) \quad (\xi\to \infty)$
required in Corollary \ref{corr0} is also checked with $H=H_\ell$.
Consequently, the same conclusions as in the previous example apply
for this process as well, in the respective regions determined by
the parameter $H_\ell$ at high frequencies $x> \omega_\ell$ alone.
The same result is also obtained for a more general process defined
by $f(\xi) = c \xi^{-2H-1}$ for $|\xi|\geq \omega$ and condition
(\ref{cond:f}) is only required elsewhere. 
Once again, such conclusions hold also in case of
continuously differentiable multiplicative and additive trends.
\end{example}

\begin{example} {\it Multifractional Brownian motion} (mBm) (see  Ayache {\it et al.}, 2005).
A mBm $X=(X_t, \, t\in
[0,1])$ is a Gaussian process defined by
\begin{equation}
X_t \ = \  B_{H(t)}(t) \ =  \ g(H(t))
\int_{\R} \frac{{\e}^{{\i}t x} - 1}{|x|^{H(t)+1/2}} W(\d
x),   \label{mBm}
\end{equation}
where $W(\d x)$ is the same as in (\ref{spec_V}), $H(t)$ is a
(continuous) function on $[0,1] $ taking values in $(0,1)$ and
finally, $g(H(t))$ is a normalization such that $\E X^2_t  = 1 $.
It is well-known that a mBm is locally asymptotically self-similar
at each point $t \in (0,1)$ having a fBm $B_{H(t)}$ as its tangent
process at $t$ (see Benassi {\it et al.}, 1997). This example is studied in more detail
in Bardet and Surgailis~(2010).

\end{example}

\begin{example} {\em Time-varying fractionally integrated processes}. Philippe {\it et al.}~(2006, 2008)
introduced two classes
of mutually inverse time-varying fractionally integrated filters
with discrete time and studied long-memory properties of the
corresponding filtered white noise processes. Surgailis~(2008)
extended these filters to continuous time and defined
``multifractional'' Gaussian processes $(X_t,\, t\ge 0) $ and $(Y_t,
t\ge 0)$ as follows
\begin{eqnarray}
X_t&=&\int_{\R} \Big\{\int_0^t \frac{1}{\Gamma (H(\tau)-.5)}(\tau -s)_+^{H(\tau) - 1.5}
{\e}^{A_-(s,\tau)} \d \tau \Big\} \d B(s), \label{X_t} \\
Y_t&=&\int_{\R} \frac{1}{\Gamma (H(s)+.5)}\Big\{ (t -s)_+^{H(s) -
.5} {\e}^{-A_+(s,t)} - (-s)_+^{H(s)-.5} {\e}^{-A_+(s,0)}  \Big\} \d
B(s), \label{Y_t}
\end{eqnarray}
where $s_+^\alpha := s^\alpha {\mathbf 1}(s >0)$, $B$ is a Brownian
motion,
$$
A_-(s,t):= \int_s^t \frac{H(u)-H(t)}{t-u}\d u, \qquad A_+(s,t):=
\int_s^t \frac{H(s)-H(v)}{v-s} \d v \qquad (s<t)
$$
and where  $H(t), t\in {\R}$ is a general function taking values in
$(0,\infty)$ and satisfying some weak additional conditions.
Surgailis~(2008) studied small and large scale limits of $(X_t)$ and
$(Y_t)$ and showed that these processes resemble a fBm with Hurst
parameter $H=H(t)$ at each point $t \in {\R}$ (i.e., admit a fBm as
a tangent process) similarly to the mBm in the previous example. The
last paper also argues that these processes present a more natural
generalization of fBm than the mBm and have nicer dependence
properties of increments. We expect that the assumptions ${\bf
(A.1)}, {\bf (A.1)'}, {\bf (A.2)}_p$ can be verified for
(\ref{X_t}), (\ref{Y_t}); however, this question requires further
work.

\end{example}

\section{Processes with independent increments} \label{Levy}

In this section, we assume that $X = (X_t,\, t\ge 0)$ is a
(homogeneous) L\'evy process, with a.s. right continuous
trajectories, $X_0 =0$.  It is well-known that if the generating triplet of $X$ satisfies certain conditions
(in particularly, if the L\'evy measure $\nu $ behaves regularly at the origin with index $\alpha \in (0,2)$), then
$X$ has a tangent process $Y$ which is $\alpha-$stable L\'evy process. A natural question is to estimate the parameter
$\alpha $ with the help of the introduced statistics $R^{p,n}$.
Unfortunately,  the limit of these statistics as defined in (1.5) through the tangent process
depends also on the skewness parameter $\beta \in [-1,1]$ of the $\alpha-$stable tangent process $Y$ and
so this limit cannot be used for determining of $\alpha $ if $\beta $ is unknown.

In order to avoid this difficulty, we shall slightly modified our ratio statistic, as follow. Observe first
that the second differences $\Delta^{2,n}_k X$ of L\'evy process have a {\it symmetric} distribution (in contrast to the first
differences $\Delta^{1,n}_k X$ which are not necessary symmetric).
For notational simplicity we shall  assume
in this section that $n$ is {\it even}.
The modified statistic
\begin{eqnarray*}
\tilde R^{2,n}&:=&\frac{1}{n/2-1}\sum_{k=0}^{(n-4)/2} \psi\big(\Delta^{2,n}_{2k} X, \Delta^{2,n}_{2k+2}X\big), \qquad \psi(x,y) := \frac{|x+y|}{|x|+|y|}
\end{eqnarray*}
is written in terms of ``disjoint'' (independent) second order  increments $(\Delta^{2,n}_{2k} X,\Delta^{2,n}_{2k+2}X)$ having
a symmetric joint distribution. Instead of extending general result of Proposition 2.2 to $\tilde R^{2,n}$, we shall directly obtain
its convergence under suitable assumptions on $X$.
Note first
\begin{equation}
\E \tilde R^{2,n} \
=\ \E \psi \big(X^{(2)}_{1/n} - X^{(1)}_{1/n}, X^{(4)}_{1/n} - X^{(3)}_{1/n}\big),
\end{equation}
where $X^{(i)}, i=1,\cdots, 4 $ are independent copies of $X$. Note that $1/2
\le \E \tilde R^{2,n} \le 1 $ since
\begin{eqnarray*}
\E \psi \big(X^{(2)}_{1/n} - X^{(1)}_{1/n}, X^{(4)}_{1/n} - X^{(3)}_{1/n}\big)&\hspace{-3mm}
\geq & \hspace{-3mm} \P(X^{(2)}_{1/n} - X^{(1)}_{1/n} \ge 0,
X^{(4)}_{1/n} - X^{(3)}_{1/n} \ge 0)\\
&& \hspace{4cm}+ \ \P(X^{(2)}_{1/n} - X^{(1)}_{1/n} <0, X^{(4)}_{1/n} - X^{(3)}_{1/n} <0) \\
&\hspace{-3mm} \geq &\hspace{-3mm}
\P^2(X^{(2)}_{1/n} - X^{(1)}_{1/n}\ge 0) + \P^2 (X^{(2)}_{1/n} - X^{(1)}_{1/n}<0) \ge 1/2 .
\end{eqnarray*}

\begin{proposition}\label{prop4.1}
Let there exists a limit
\begin{equation}
\lim_{n\to \infty} \E \tilde R^{2,n} =  \tilde \Lambda.  \label{R_mean_conv}
\end{equation}
Then
\begin{equation}
\tilde R^{2,n} \limitepsn \tilde \Lambda. \label{as}
\end{equation}
\end{proposition}
\noindent {\it Proof.} Write $\tilde R^{2,n} = \E \tilde R^{2,n} + (n/2-1)^{-1}Q_n,
$ where $Q_n$ is a
sum of centered 1-dependent  r.v.'s which are bounded by 1 in absolute value. Therefore
$\E ((n/2-1)^{-1} Q_n)^4 = O(n^{-2}) $ and
the a.s. convergence $(n/2-1)^{-1}Q_n \to 0 $ follows by
the Chebyshev inequality. \hfill  $\Box$

\bigskip

Next we discuss conditions on $X$ for the convergence in
(\ref{R_mean_conv}). Recall that the distribution of $X_t$ is
infinitely divisible and its characteristic function is given by
\begin{equation}
\E \e^{{\i}\theta X_t} = \exp \Big\{t \Big({\i}\gamma \theta - \frac{1}{2}a^2 \theta^2 + \int_{\R}
\big({\e}^{{\i}u \theta } -1-  {\i}u\theta \1(|u|\le 1)\big) \nu ({\d}u)
\Big)\Big\},  \quad \theta \in {\R}, \label{ch_f}
\end{equation}
where $\gamma \in {\R}, a \ge 0 $ and $\nu $ is a measure on $\R$
such that $\int_{\R} \min (u^2, 1) \nu({\d}u) < \infty$. The
triplet $(a, \gamma, \nu)$ is called the generating triplet of $X$
(Sato~(1999)). Let $X^{(i)}, i=1,2 $ be independent copies of $X$. Note $W_t := X^{(1)}_t - X^{(2)}_t $ is a L\'evy process having the characteristic function
\begin{equation}
\E \e^{{\i}\theta W_t} = \exp \Big\{t \Big(-a^2 \theta^2 + 2\int_0^\infty
{\rm Re}\big(1- {\e}^{{\i}u \theta }\big) {\d}K(u)
\Big)\Big\},  \quad \theta \in {\R}, \label{chY}
\end{equation}
where
$$
K(u)\ := \ \nu ((-\infty,-u]\cup [u, \infty))
$$
is monotone nonincreasing on $(0,\infty)$.
Introduce the following condition: there exist $0<
\alpha  \leq  2$ and $c>0 $
such that
\begin{equation}
K(u) \sim \frac{c}{u^\alpha}, \quad u \downarrow 0. \label{K}
\end{equation}
It is clear that if such number $\alpha$ exists then
$\alpha := \inf\big \{ r \geq 0 : \int_{ |x| \leq 1} |x|^r\nu (dx) < \infty \big \}$
is the so-called  fractional order or the Blumenthal-Getoor index of the Lévy process $X$.\\
Let $Z_\alpha $ be a standard $\alpha-$stable r.v. with characteristic function
$
\E \e^{{\i}\theta Z_\alpha} = \e^{- |\theta|^\alpha}
$
and $Z^{(i)}_{\alpha}, i=1,2,3$ be independent copies of $Z_\alpha $.

\medskip
\begin{proposition}
\label{prop4.2}
Assume either $a>0$ or else, $a=0$ and condition (\ref{K}) with $0< \alpha \le 2$ and $c>0$. Then $t^{-1/\alpha} (X^{(1)}_t - X^{(2)}_t) \limiteloit \tilde c \, Z_\alpha $  with $\tilde c$ depending on $c$, and
(\ref{R_mean_conv}), (\ref{as}) hold with
$$
\tilde \Lambda \equiv \tilde \Lambda(\alpha) := \E \psi \big(Z^{(1)}_\alpha, Z^{(2)}_\alpha \big).
$$
Moreover, with $\displaystyle
\tilde \sigma^2(\alpha) := 2{\rm var}\big(\psi \big(Z^{(1)}_\alpha, Z^{(2)}_\alpha \big)\big) + 4 {\rm cov}
\big(\psi \big(Z^{(1)}_\alpha, Z^{(2)}_\alpha \big),  \psi\big(Z^{(2)}_\alpha, Z^{(3)}_\alpha\big)\big)$,
\begin{equation} \label{tildeCRT}
\sqrt{n}(\tilde R^{2,n}- \E \tilde R^{2,n}) \limiteloin {\cal N}(0, \tilde \sigma^2(\alpha)).
\end{equation}

\end{proposition}

\noindent {\it Proof.} Relation $t^{-1/\alpha} W_t = t^{-1/\alpha} (X^{(1)}_t - X^{(2)}_t) \limiteloit \tilde c \, Z_\alpha $
is an easy consequence of the assumptions of the proposition and the general criterion of weak convergence
of infinitely divisible distributions in  Sato~(1991, Theorem 8.7). It implies (\ref{R_mean_conv}) by the fact that  $\psi$
is a.e. continuous on ${\R}^2$.  Since $\tilde R^{2,n}$ is a sum of 1-dependent stationary and bounded r.v.'s, the central limit theorem
in (\ref{tildeCRT}) follows from convergence of the variance:
\begin{equation} \label{tildesigma}
n \, \Var (\tilde R^{2,n}) \to \tilde \sigma^2(\alpha);
\end{equation}
see e.g. Berk~(1973). Rewrite $\tilde R^{2,n} = (n/2-1)^{-1}\sum_{k=0}^{(n-4)/2} \tilde\eta_{n}(k), \ \tilde \eta_{n}(k):= \psi\big(\Delta^{2,n}_{2k} X, \Delta^{2,n}_{2k+2}X\big)$.
We have
$$
n \, \Var(\tilde R^{2,n})= \frac{n}{n/2-1} \, \Var(\tilde \eta_n(0)) + \frac{2n(n/2-2)}{(n/2-1)^2} \Cov(\tilde \eta_n(0), \tilde \eta_n(1)),
$$
where $\Var(\tilde \eta_n(0))\limiten \Var \big(\psi \big(Z^{(1)}_\alpha, Z^{(2)}_\alpha \big)\big), \
\Cov (\tilde \eta_n(0), \tilde \eta_n(1)) \limiten \Cov \big(\psi \big(Z^{(1)}_\alpha, Z^{(2)}_\alpha \big),  \psi\big(Z^{(2)}_\alpha, Z^{(3)}_\alpha\big)\big)$ similarly as in
the proof of (\ref{R_mean_conv}) above. This proves (\ref{tildesigma}) and the proposition. \hfill $\Box$
\begin{center}
\begin{figure}[h]
\begin{center}
\includegraphics[width=7 cm,height=3.6cm]{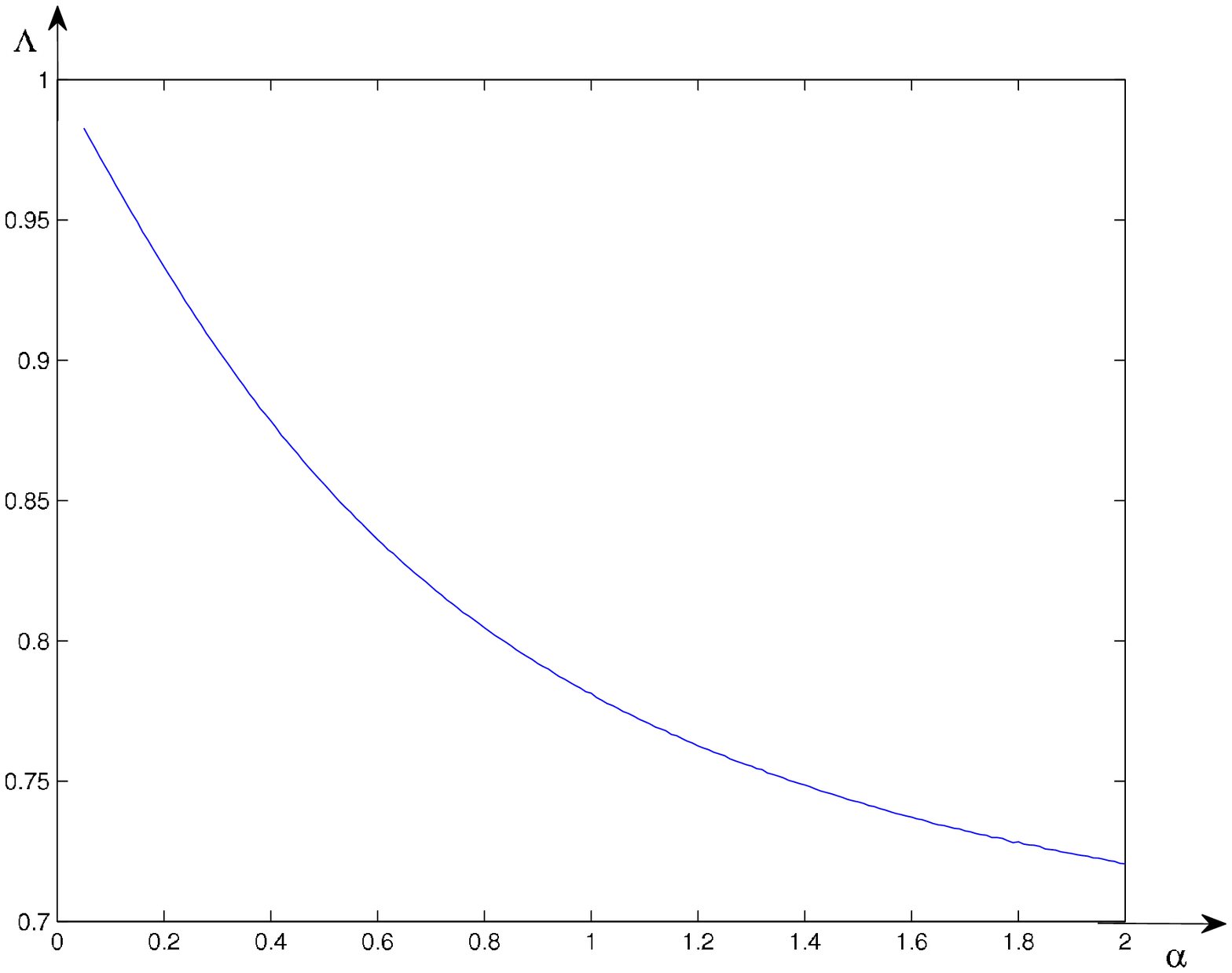}\hspace{2mm}
\includegraphics[width=8 cm,height=4cm]{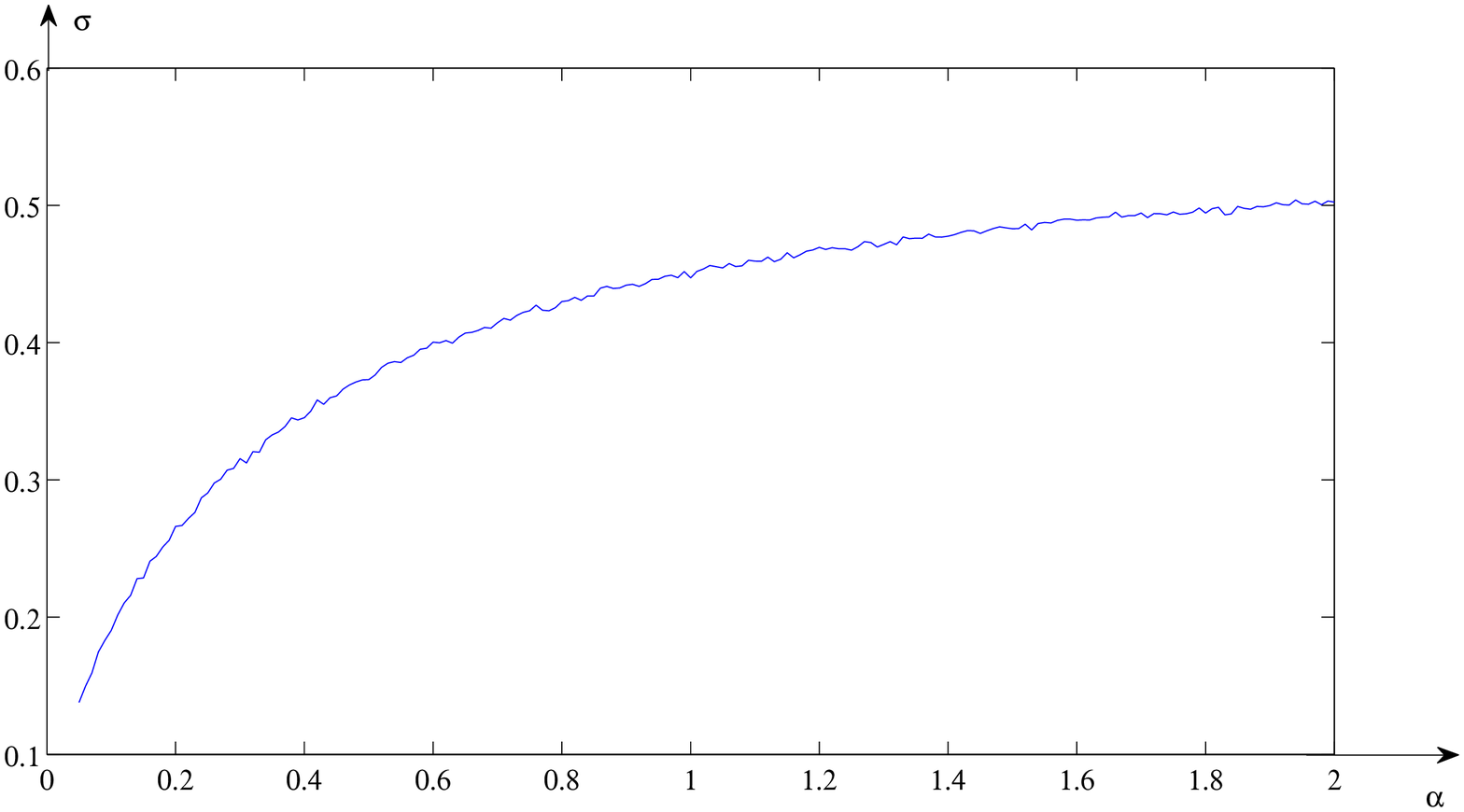}
\caption{The graphs of $\alpha \mapsto \tilde \Lambda(\alpha) =
\E \frac{|Z^{(1)}_{\alpha}+ Z^{(2)}_{\alpha}|}{|Z^{(1)}_{\alpha}| + |Z^{(2)}_{\alpha}|}$ (left) and $\alpha \mapsto \tilde \sigma(\alpha)$ (right) for a process with independent increments. } \label{Figure2}
\end{center}
\end{figure}
\end{center}
The graph of $\tilde\Lambda (\alpha)$ 
is given in Figure \ref{Figure2}. Note that $\tilde\Lambda (2)= \Lambda_1(1/2)\simeq 0.72$: this is the case of Brownian motion. \\
In order to evaluate the decay rate of the bias $\E \tilde R^{2,n} -  \tilde \Lambda(\alpha)$ we need a uniform convergence rate in
Lemma \ref{lemmaFG}, below, for
$$
\|F_n - G_\alpha\|_\infty := \sup_{x\in \R} |F_n(x) - G_\alpha(x)|, \quad
F_n(x) := \P(n^{1/\alpha} W_{1/n}\le x), \quad G_\alpha(x) := \P(\tilde Z_\alpha \le x),
$$
where $\tilde Z_\alpha := \tilde c Z_\alpha $ is the limiting $\alpha-$stable r.v. in Proposition \ref{prop4.2} and
$(W_t, t\ge 0) $ is the symmetric Levy process with characteristic function as in (\ref{chY}).  The proof of Lemma \ref{lemmaFG}
is given in Annexe.

\begin{lemma} \label{lemmaFG}
(i) \, Let $a=0$ and $K$ satisfy (\ref{K}). Denote $K_1(u):= K(u) - cu^{-\alpha}, $    $|K_1|(u) := \int_u^\infty |{\d}K_1(v)| $, the variation of $K_1$ on $[u, \infty)$.
Moreover, assume that there exist some constants $\beta, \delta >0$ such that
\begin{equation} \label{KK}
|K_1|(u) = O(u^{- (\alpha -\beta)_+}) \quad (u \to 0), \quad |K_1|(u) = O(u^{-\delta}) \quad (u \to \infty),
\end{equation}
where $x_+ := \max(0,x)$. Then
\begin{eqnarray}\label{FG}
\|F_n-G_\alpha\|_\infty&=&\cases{O(n^{-\beta/\alpha}), &if $\beta < \alpha$, \cr
O(n^{-1} \log n), &if $\beta = \alpha$, \cr
O(n^{-1}), &if $\beta > \alpha $. \cr}
\end{eqnarray}

\noindent (ii)\, Let $a >0$ and $K$ satisfy
\begin{equation} \label{KKK}
K(u) = O(u^{- \alpha}) \quad (u \to 0), \quad K(u) = O(u^{-\delta}) \quad (u \to \infty),
\end{equation}
for some $0\le \alpha < 2, \ \delta >0$. Then
\begin{eqnarray}\label{FG2}
\|F_n-G_\alpha\|_\infty&=&\cases{O(n^{-1+\alpha/2}), &if $\alpha>0$, \cr
O(n^{-1} \log n), &if $\alpha=0$. \cr}
\end{eqnarray}

\end{lemma}

\medskip

\begin{proposition} \label{RBias} Assume either $a>0$ or else $a=0$ and condition (\ref{K}). Then for any $\alpha \in (0,2]$
\begin{equation} \label{bias}
|\E \tilde R^{2,n} -  \tilde \Lambda(\alpha)| \ \le \ 2C \|F_n - G_\alpha\|_\infty, \qquad C :=  \int_0^\infty (1+z)^{-2} {\d}z < \infty.
\end{equation}
\end{proposition}

\noindent {\it Proof.} Let $\tilde \psi(x,y) := |x-y|/(x+y), \  x,y >0$, and let $F_n, G_\alpha$ be the same as in Lemma \ref{lemmaFG}.
Similarly as in Vai\v ciulis (2009, proof of Th.~1), write
\begin{eqnarray*}
\E \tilde R^{2,n} -  \tilde \Lambda(\alpha)&=&2\int_0^\infty \int_0^\infty \tilde \psi(x,y) ({\d}F_n(x) {\d} F_n(y) - {\d}G_\alpha(x) {\d}G_\alpha(y)) \ = \ 2(W_1 + W_2),
\end{eqnarray*}
where $W_1:=\int_0^\infty \int_0^\infty \tilde \psi(x,y) {\d}F_n(x)({\d} F_n(y) - {\d}G_\alpha(y)), \
W_2:=\int_0^\infty \int_0^\infty \tilde \psi(x,y) {\d}G_\alpha(y)( {\d} F_n(x) - {\d}G_\alpha(x)). $  Integrating by parts yields
\begin{eqnarray*}
|W_1|&=&2\int_0^\infty |x| {\d}F_n(x) \int_0^\infty |F_n(y) - G_\alpha (y)| \frac{{\d}y}{(x+y)^2} \\
&\le&
2 \|F_n- G_\alpha\|_\infty \int_0^\infty |x| {\d}F_n(x) \int_0^\infty  \frac{{\d}y}{(x+y)^2} \
= \  C \|F_n- G_\alpha\|_\infty
\end{eqnarray*}
since $\int_0^\infty {\d}F_n(x) = 1/2 $.
A similar estimate holds for $W_2$. This proves  (\ref{bias}). \hfill $\Box$

\medskip

Propositions \ref{prop4.2}, \ref{RBias}, and Lemma \ref{lemmaFG}, together with  the Delta-method, yield the following corollary.

\begin{corollary}\label{estimaalpha} Let $a$ and $K$ satisfy either the assumptions of Lemma 5.1(i) with $\beta > \alpha/2$, or the assumptions
of Lemma 5.1(ii). Then
$$
\sqrt{n}(\tilde R^{2,n} - \tilde \Lambda(\alpha))\limiteloin {\cal N}(0, \tilde \sigma^2(\alpha)).
$$
Moreover, if we define $\widehat \alpha_n:=\tilde \Lambda^{-1}(\tilde R^{2,n})$, then
$$
\sqrt{n}\big(\widehat \alpha_n - \alpha\big)\limiteloin {\cal N}(0, \tilde s^2(\alpha)),
$$
where $\tilde s^2(\alpha ) := \big[\frac{\partial  {\tilde \Lambda}}{\partial \alpha} (\alpha)\big]^{-2} \tilde \sigma^2(\alpha), \  0 < \alpha \le 2 $.

\end{corollary}

There exist very few papers concerning  estimation of $\alpha$ in such a semiparametric frame. Nonparametric estimation of parameters of
Lévy processes based on the empirical characteristic function has recently been considered in Neumann and Reiß~(2009) and
Gugushvili~(2008), but the convergence rates there are $(\log n)^\kappa$ with $\kappa>0$. Ait Sahalia and Jacod~(2009) have proposed an estimator of the degree of activity
of jumps (which is identical to the fractional order in the case of a Lévy process) in a general semimartingale framework using small increments of high frequency data. However from the generality of their model, the convergence rate of the estimator is not rate efficient (in fact smaller than $n^{1/5}$). A recent paper of Belomestny~(2010) provides an efficient data-driven procedure to estimate $\alpha$ using a spectral approach but in a different semiparametric frame from ours. Thus,
Corollary \ref{estimaalpha} appears as a new and interesting result since the estimator $\widehat \alpha_n$ follows a $\sqrt n$-central limit theorem.

\section{Annexe : proofs} \label{proofs}

\subsection*{\it Proof of Lemma \ref{Lem32}}

Let $\delta^2 := \max_{i=1,2} \E \xi_i^2$. If $\delta^2 \ge
1/2$, then (\ref{Z_bdd}) holds since the l.h.s. of
(\ref{Z_bdd}) does not exceed $1$. Let $\delta^2< 1/2$ in the sequel.
Write $U:= \psi (Z_1+\xi_1, Z_2+ \xi_2) - \psi(Z_1,Z_2) = U_\delta +
U^c_\delta$,
\begin{eqnarray*}
U_\delta&:=&U {\mathbf 1}(A_\delta) \ = \  (\psi (Z_1+\xi_1, Z_2+ \xi_2) - \psi(Z_1,Z_2)) {\mathbf 1}(A_\delta) ,\\
U^c_\delta&:=&U {\mathbf 1}(A^c_\delta) \ = \ (\psi (Z_1+\xi_1, Z_2+
\xi_2) - \psi(Z_1,Z_2)) {\mathbf 1}(A^c_\delta),
\end{eqnarray*}
where ${\mathbf 1}(A_\delta) $ is the indicator of the event
$$
A_\delta := \{|Z_1| > \delta^{2/3},  |Z_2|> \delta^{2/3},  |\xi_1| <
\delta^{2/3}/2, |\xi_2|< \delta^{2/3}/2\},
$$
and ${\mathbf 1}(A^c_\delta) = 1 - {\mathbf 1}(A_\delta) $ is the
indicator of the complementary event $A^c$.  Clearly,
\begin{eqnarray*}
\E |U^c_\delta| &\le& 2(\P(|Z_1|< \delta^{2/3}) + \P(|Z_2|<\delta^{2/3}) +
\P(|\xi_1|\ge \delta^{2/3}/2) + \P(|\xi_2|\ge \delta^{2/3}/2)\\
&\le&\frac{4}{\sqrt{2\pi}}\, \delta^{2/3} + 2\, \frac{\E
\xi^2}{\delta^{4/3}}\ \le 4\ \delta^{2/3}.
\end{eqnarray*}
It remains to estimate $\E |U_\delta|$. By the mean value theorem,
$$
|U_\delta| \le \left(|\xi_1| \sup_{D}|\psi_{x_1}| + |\xi_2| \sup_{D}
|\psi_{x_2}|\right) {\mathbf 1}(A_\delta),
$$
where
\begin{eqnarray*}
\sup_D |\psi_{x_1}| &:=& \sup \{|\partial \psi (x_1,x_2)/\partial x_1|: |x_i-Z_i| \le |\xi_i|,  i=1,2 \}, \\
\sup_D |\psi_{x_2}| &:=& \sup \{|\partial \psi (x_1,x_2)/\partial
x_2|:  |x_i- Z_i| \le |\xi_i|,  i=1,2\}.
\end{eqnarray*}
Therefore
\begin{eqnarray*}
\E |U_\delta| &\le& \E^{1/2}(\xi_1^2) \E^{1/2}\big [
(\sup_{D}|\psi_{x_1}|)^2
{\mathbf 1}(A_\delta)\big ] + \E^{1/2}  (\xi_2^2) \E^{1/2} \big [ (\sup_{D}|\psi_{x_2}|)^2 {\mathbf 1}(A_\delta)\big ] \\
&\le& \delta \left(\E^{1/2} \big [ (\sup_{D}|\psi_{x_1}|)^2 {\mathbf
1}(A_\delta) \big ]+ \E^{1/2}\big [ (\sup_{D}|\psi_{x_2}|)^2 {\mathbf
1}(A_\delta)\big ] \right).
\end{eqnarray*}
Next,
$$
|\psi_{x_i} (x_1,x_2)| = \frac{|{\rm sgn}(x_1+x_2)(|x_1|+|x_2|) -
(|x_1+x_2|){\rm sgn}(x_i)|} {(|x_1|+|x_2|)^2 } \  \le \
\frac{2}{|x_1|+|x_2|}.
$$
Therefore
\begin{eqnarray*}
\sup_{D}|\psi_{x_i}|^2 {\mathbf 1}(A_\delta)
&\le&4 \, \sup \big \{ (|x_1|+|x_2|)^{-2}:  |x_i -Z_i| \le \delta^{2/3}/2, |Z_i| > \delta^{2/3}, i=1,2 \big \} \\
&\le&16 \,  (|Z_1|+|Z_2|)^{-2} {\mathbf 1}(|Z_i| > \delta^{2/3}, i=1,2),
\end{eqnarray*}
implying
$$
\E \big [(\sup_{D}|\psi_{x_i}|)^2 {\mathbf 1}(A_\delta)\big ]  \le 16 \, \E \left[
\frac{1}{(|Z_1|+|Z_2|)^2}; |Z_i| > \delta^{2/3}, i=1,2 \right] \ \le
\ C(\delta),
$$
where
\begin{eqnarray*}
C(\delta) &\le& \frac{16}{2\pi} \int_{\{x_1^2+x_2^2 >
\delta^{4/3}\}}
\frac{1}{x_1^2 + x_2^2} {\e}^{-(x_1^2+x_2^2)/2} {\d}x_1 {\d}x_2 \\
&\leq &16\int_{\delta^{4/3}}^\infty r^{-1} {\e}^{-r^2/2} {\d}r \  \le \
16\, (1 +(4/3)|\log \delta|).
\end{eqnarray*}
Hence $\E |U_\delta| \le 8 \, \delta (1+ (4/3) \, |\log \delta|)^{1/2}, ~~ \E
|U^c_\delta| \le 4\, \delta^{2/3}$. It remains to use $x(1+ (4/3) \,|\log
x|)^{1/2} \le 2x^{2/3}$ for all $0<x \le 1$. \hfill $\Box$.

\subsection*{\it Proof of Property \ref{prop31}}

We use the following identity: for any reals $x_1, \cdots, x_j$,
$$
x_1 x_j = \frac{1}{2}\bigg\{\Big(\sum_{k=1}^j x_k\Big)^2 + \Big(\sum_{k=2}^{j-1} x_k \Big)^2
- \Big(\sum_{k=1}^{j-1} x_k \Big)^2 - \Big(\sum_{k=2}^j x_k \Big)^2 \bigg\}.  \label{xx}
$$
In particular,
\begin{eqnarray*}
\E \left[\Delta^{1,n}_{[nt]}X \Delta^{1,n}_{j+[nt]}X \right]
&=&\frac{1}{2}\bigg\{\E \big(X_{j+1+ [nt]} - X_{[nt]}\big)^2 +  \E \big(X_{j-1+ [nt_*]} - X_{[nt_*]}\big)^2\\
&-&\E \big(X_{j+ [nt]} - X_{[nt]}\big)^2 - \E \big(X_{j+ [nt_*]} - X_{[nt_*]}\big)^2\bigg\},
\end{eqnarray*}
where $t_* := t+ (1/n)$ (so that $[nt_*] = [nt] + 1)$. Then, using
{\bf (A.1)} and the notation $u_n$ for a sequence tending to $0$
as $n \to \infty $ uniformly in $t$ and all $|j| < J$, where $J$
is a fixed number, we obtain
\begin{eqnarray*}
\E \left[\Delta^{1,n}_{[nt]}X \Delta^{1,n}_{j+[nt]}X \right]
&=&\frac{1}{2}\bigg\{c(t)\Big(\mbox{$\frac{j+1}{n}$}\Big)^{2H(t)}\left(1+u_n\right)
+ c(t_*)\Big(\mbox{$\frac{j-1}{n}$}\Big)^{2H(t_*)}\left(1+u_n\right) \\
&&\hskip.5cm - \
c(t)\Big(\mbox{$\frac{j}{n}$}\Big)^{2H(t)}\left(1+u_n\right)
- c(t_*)\Big(\mbox{$\frac{j}{n}$}\Big)^{2H(t_*)}\left(1+u_n\right) \bigg\}  \\
&=&\frac{c(t)}{2} \Big\{ \Big(\mbox{$\frac{j+1}{n}$}\Big)^{2H(t)} +
 \Big(\mbox{$\frac{j-1}{n}$}\Big)^{2H(t)}
-2\Big(\mbox{$\frac{j+1}{n}$}\Big)^{2H(t)}
\Big\}\left(1+u_n\right)
\end{eqnarray*}
since $c(t_*) - c(t) = u_n $ and $\Big(\frac{j}{n}\Big)^{2H(t_*)}
= \Big(\frac{j}{n}\Big)^{2H(t)}\left(1+u_n\right) $ follows from
(\ref{H_cont}). This proves (\ref{R_j1}) for $p=1$. Relation
(\ref{R_j1}) for $p=2$ follows analogously. Relation (\ref{R_j11})
also follows by the same argument  and the fact that $c(t_*) -
c(t) = u_n/\sqrt{n} $ and $\Big(\frac{j}{n}\Big)^{2H(t_*)} =
\Big(\frac{j}{n}\Big)^{2H(t)}\left(1+u_n/\sqrt{n}\right) $ hold in
view of Assumption {\bf (A.1)'}. Property \ref{prop31} is proved.
\hfill $\Box$

\subsection*{\it Proof of Property \ref{prop_V}}

With condition $V(t) \sim ct^{2H} \  (t \to 0)$ in mind,
inequality (\ref{rho_dom2}) reduces to
\begin{eqnarray}
\left|
V\left(\mbox{$\frac{k+1}{n}$}\right)+ V\left(\mbox{$\frac{k-1}{n}$}\right) - 2V\left(\mbox{$\frac{k}{n}$}\right)\right|
&\le &
Cn^{-2H+\theta} k^{-\gamma} \quad (p=1, \ 2\le k \le n), \label{V_1}\\
\left|
V\left(\mbox{$\frac{k+2}{n}$}\right) - 4V\left(\mbox{$\frac{k+1}{n}$}\right)+ 6V\left(\mbox{$\frac{k}{n}$}\right) -
4V\left(\mbox{$\frac{k-1}{n}$}\right) + V\left(\mbox{$\frac{k-2}{n}$}\right)\right|
&\le &
Cn^{-2H+\theta} k^{-\gamma} \quad (p=2, \ 4\le k \le n). \label{V_2}
\end{eqnarray}
The left hand side of (\ref{V_1}) can be written and estimated as
\begin{eqnarray*}
\left|\int_0^{1/n}  \int_0^{1/n} V''(t-s+(k/n)) \d t \d s\right|
&\le&C\int_0^{1/n} \int_0^{1/n} |t-s + (k/n)|^{-\kappa} \d t \d s \\
&\leq &Cn^{\kappa-2} \int_0^1 \int_0^1 |t-s+k|^{-\kappa} \d t \d s \\
&\le&Cn^{\kappa-2} k^{-\kappa} \  =  \  Cn^{\theta -2H} k^{-\gamma}
\end{eqnarray*}
where $\gamma = \kappa> 1/2 $ and $\theta = \kappa +2H-2 \in [0, \gamma/2) $ since $\kappa < 4-4H$.
This proves part (i).  Part (ii) follows similarly, by writing
the left hand side of (\ref{V_2}) as
\begin{eqnarray*}
\left|\int_0^{1/n}\hskip-.3cm  \cdots \int_0^{1/n} V^{(4)}(t-s+ u-  v +(k/n)) \d t \d s \d u \d v\right|
&\le&C\int_0^{1/n}\hskip-.3cm \cdots \int_0^{1/n} |t-s + u - v + (k/n)|^{-\kappa} \d t \d s  \d u \d v \\
&\le&Cn^{\kappa-4} k^{-\kappa} \ = \ Cn^{\theta-2H} k^{-\gamma},
\end{eqnarray*}
where $\gamma = \kappa > 1 $ and $\theta = \kappa +2H- 4 \in [0,
\kappa/2)$ since $\kappa < 8-4H$. Property \ref{prop_V} is proved.
\hfill $\Box$

\subsection*{\it Proof of Property \ref{prop_f}}

{}From (\ref{spec_V}) we have
\begin{eqnarray*}
\E \big [\Delta^{p,n}_0 X \  \Delta^{p,n}_{j} X \big ]
&=&\int_{\R} \left|\e^{\i (x/n)} -1\right|^{2p} \e^{\i x(j/n)} f(x) \d x \ = \
2^{1+p}\int_0^\infty (1- \cos (x/n))^p \cos(xj/n) f(x) \d x \nonumber \\
&=&2^{1+p}(I_1 + I_2),
\end{eqnarray*}
with
\begin{eqnarray}
|I_1|&=&\Big |\int_0^K (1- \cos (x/n))^p \cos(xj/n) f(x) \d x \Big | \nonumber  \\
&=&C \Big | \int_0^{K} (x/n)^{2p} f(x) \d x\Big |  \ = \ \Big | C
n^{-2p} \int_0^K x^{2p} f(x) \d x \Big | \ \le \  C n^{-2p};
\label{I1} \\
\mbox{and}~~~I_2&=&n \int_{K/n}^\infty (1- \cos (x))^p \cos(xj) f(nx) \d x \nonumber \\
&=&(n/j)\int_{K/n}^\infty (1- \cos (x))^p f(nx) \d \sin(xj) \nonumber \\
&=&-(n/j) \int_{K/n}^\infty \sin (xj) \left((1- \cos (x))^p
f(nx)\right)^\prime_x \d x +  O(n^{-2p}) \label{fourier}
\end{eqnarray}
since $ |(n/j) f(K) (1- \cos(K/n))^p \sin(jK/n)| \  \le \
C(n/j)(K/n)^{2p} |jK/n| \le Cn^{-2p} $ for any $K>0$ fixed. \\
~\\
Let $p=1$. The last integral can be rewritten as
$\int_{K/n}^\infty \sin (xj) \left((1- \cos (x))
f(nx)\right)^\prime_x \d x = \tilde I_1 + \tilde I_2, $ where
\begin{eqnarray*}
|\tilde I_1| &\le&\int_{K/n}^\infty |\sin (xj) \sin (x) f(nx)| \d
x \  \le \  C\int_{K/n}^\infty |\sin(xj) \sin(x)| (nx)^{-2H-1} \d
x \  \le \
 \sum_{q=1}^3 I_{1q},
\end{eqnarray*}
where we used the fact that $f(x) \le Cx^{-2H-1} \  (x \ge K) $;
see condition (\ref{spec_f}), and where
\begin{eqnarray*}
|I_{11}|&=&\int_0^{1/j} |\sin (xj) \sin (x)| (nx)^{-2H-1} \d x  \   \le \
Cj \int_0^{1/j} x^2 (nx)^{-1-2H} \d x \  \le \  Cj^{2H-1}n^{-1-2H}, \\
|I_{12}|&=&\int_{1/j}^1  |\sin (xj) \sin (x)| (nx)^{-2H-1} \d x  \  \le \
C\int_{1/j}^1x (nx)^{-1-2H} \d x \  \le \  Cj^{2H-1} n^{-1-2H}, \\
|I_{13}|&=&\int_{1}^\infty  |\sin (xj) \sin (x)|(nx)^{-2H-1} \d x
\  \le \  C\int_{1}^\infty (nx)^{-1-2H} \d x \  \le \  Cn^{-1-2H}.
\end{eqnarray*}
Similarly,  using (\ref{spec_f}),
\begin{eqnarray*}
|\tilde I_2|&\le&n \int_{K/n}^\infty |\sin (xj) (1- \cos (x))
f'(nx)| \d x  \  \le  \   C\int_{K/n}^\infty |\sin(xj)(1- \cos(x)|
(nx)^{-2H-2} \d x \  \le \ \sum_{q=1}^3 I_{2q},
\end{eqnarray*}
where
\begin{eqnarray*}
 I_{21} &=& n\int_0^{1/j} |\sin (xj) (1- \cos (x))| (nx)^{-2H-2}
\d x
\ \le \   C n j \int_0^{1/j} x^3 (nx)^{-2-2H} \d x  \  \le \  C  j^{2H-1}  n^{-1-2H}, \\
I_{22} &=&  n\int_{1/j}^1 |\sin (xj) (1- \cos (x))| (nx)^{-2H-2}
\d x  \ \le \  Cn \int_{1/j}^1 x^2 (nx)^{-2-2H} \d x  \
\le \  C  j^{2H-1}  n^{-1-2H}, \\
I_{23} &=&n\int_{1}^\infty |\sin (xj) (1- \cos (x))| (nx)^{-2H-2}
\d x  \ \le \  Cn \int_{1}^\infty (nx)^{-2-2H} \d x \ \le \ C
n^{-1-2H}.
\end{eqnarray*}
We finally obtain, for $1\le j \le n$,
\begin{eqnarray*}
\left|\E \big [\Delta^{1,n}_0 X \  \Delta^{1,n}_{j} X \big ] \right|
&\le&C(n/j) j^{2H-1} n^{-1-2H} + O(n^{-2}) \  \le \ Cn^{-2H} j^{2H-2}
\end{eqnarray*}
implying (for $0<H<3/4$)  {\bf (A.2)}$_1$ with $\theta = 0$ and
$\kappa = 2-2H> 1/2 $.  For $p=2$, the estimation of the integral
in (\ref{fourier}) is completely similar, resulting in the bound
\begin{eqnarray*}
\left|\E \big [\Delta^{2,n}_0 X \  \Delta^{2,n}_{j} X \big ] \right|
&\le&C(n/j) n^{-1-2H} + O(n^{-4}) \  \le  \ Cn^{-2H} j^{-1}
\end{eqnarray*}
for any $0<H<1$, or   {\bf (A.2)}$_2$ with $\theta = 0$ and $\kappa = 1> 1/2 $.  \hfill
$\Box$

\subsection*{\it Sketch of the proof of Theorem \ref{psgauss}}

The proof of Theorem \ref{psgauss} is based on the moment inequality in Lemma \ref{lemgauss}, below,  which
extends a similar inequality in (Taqqu, 1977, Lemma 4.5) to vector-valued nonstationary Gaussian
processes. The proof of Lemma \ref{lemgauss} uses the diagram formula and is given in Bardet and Surgailis (2009).
To formulate this lemma, we need  the following definitions.  Let ${\mbf X}$ be a standard Gaussian vector in ${\R}^\nu \, (\nu \ge 1)$ and
let $\L^2({\mbf X})$ denote the Hilbert space  of measurable functions
$f: {\R}^\nu \to {\R}$ satisfying $\|f\|^2 :=  \E (f({\mbf
X}))^2 < \infty$. 
Let $\L^2_0({\mbf X}) = \{f
\in \L^2({\mbf X}): \E f({\mbf X}) = 0\}$. Let $({\mbf X}_1, \cdots, {\mbf X}_N)$ be a collection of
standardized Gaussian vectors ${\mbf X}_t = (X^{(1)}_{t}, \cdots,
X^{(\nu)}_{t}) \in {\R}^\nu $ having a joint Gaussian distribution
in ${\R}^{\nu N}$. Let $\varepsilon \in [0,1]$ be a fixed number.
Following Taqqu~(1977), we call $({\mbf X}_1, \cdots, {\mbf X}_N)$
{\it $\varepsilon-$standard} if $|\E X^{(u)}_{t} X^{(v)}_{s}| \le
\varepsilon $ for any $t \neq s, 1\le t,s\le N$ and any $1\le u, v
\le \nu$. Finally, $\sum'$ denotes the sum over all distinct integers $1 \le t_1, \cdots, t_p \le N, \ t_i \ne t_j \ (i \ne j)$.

\smallskip

\begin{lem}\label{lemgauss}
Let $({\mbf X}_1, \cdots, {\mbf X}_N)$ be $\varepsilon-$standard
Gaussian vector, ${\mbf X}_t = (X^{(1)}_{t}, \cdots,X^{(\nu)}_{t})
\in {\R}^\nu \ (\nu \ge 1)$, and let $G_{j,t,N} \in \L^2({\mbf X}),
1\le j \le p \  (p\ge 2), 1\le t \le N $. For
given integers $m, N \ge 1$, define
\begin{eqnarray}
Q_N &:=&\max_{1\le t \le N} \sum_{1\le s\le N, s\neq t} \max_{1\le
u,v \le \nu} |\E X^{(u)}_{t} X^{(v)}_{s}|^m. \label{Q_N}
\end{eqnarray}
Assume that for some integer $ 0\le \alpha \le p $,
the functions $G_{1,t,N},\cdots, G_{\alpha,t,N} $ have
a Hermite rank at least equal to $m$ for any $N\ge 1, 1\le t \le
N$, and that
$\varepsilon < \frac{1}{\nu p-1}$. Then
\begin{eqnarray*}
\sum\nolimits^\prime \E|G_{1,t_1,N}({\mbf X}_{t_1}) \cdots
G_{p,t_p,N}({\mbf X}_{t_p})| &\le& C( \varepsilon, p, m, \alpha,
\nu ) K N^{p- \frac{\alpha}{2}} Q_N^{\frac{\alpha}{2} },
\end{eqnarray*}
where the constant $C( \varepsilon, p, m, \alpha, \nu ) $ depends
on $\varepsilon, p, m, \alpha, \nu $ only, and
$K := \prod_{j=1}^p \max_{1\le t \le N}  \|G_{j,t,N}\|$.
\end{lem}
~\\

\noindent {\it Sketch of the proof of Theorem \ref{psgauss}.} The convergence $\lim_{n\to \infty} \E R^{p,n} = \int_0^1
\Lambda_p(H(t)) \d t $ is easy (see the proof of Proposition 2.2).
Hence (\ref{R_n1}) follows from
\begin{equation}
\tilde R^{p,n} := R^{p,n}  - \E R^{p,n}  \limitepsn  0. \label{tilde_R1}
\end{equation}
Relation (\ref{tilde_R1}) follows from the Chebyshev Inequality
and the following bound: there exist $C, \kappa > 1$ such that for
any $n\ge 1$
\begin{equation}
\E \Big(\tilde R^{p,n}\Big)^4 \  \le \  Cn^{-\kappa}. \label{4mom1}
\end{equation}
By definition, $\tilde R^{p,n} = \frac 1 {n-p} \sum_{k=0}^{n-p-1} \tilde \eta_n(k), $ where $  \tilde \eta_n(k):=\eta_n(k) - \E \eta_n(k) $ and
$ \eta_n(k) := \   \psi(\Delta^{p,n}_{k}X, \Delta^{p,n}_{k+1}X), $  $ \psi(x,y) = |x+y|/(|x|+|y|) $ are nonlinear functions  of Gaussian
vectors  $(\Delta^{p,n}_{k}X, \Delta^{p,n}_{k+1}X)\in \R^2$  having the Hermite rank 2; however, these vectors are not $\varepsilon-$standard and therefore Lemma 1 cannot be directly applied
to estimate the l.h.s. of (\ref{4mom1}) (with $p=1, \cdots, 4, \, \nu =2$).  To this end, we first need to ``decimate'' the sum $ \tilde R^{p,n}$, as follows. (A similar
``trick'' was used in  Cs\"org\H o and Mielniczuk (1996).)
Let $\ell = [n^{\theta/\gamma}]$ be the sequence of integers increasing to $\infty$ (at a rate  $o(n^{1/2})$ by condition $\theta <
\gamma/2$) and write
\begin{eqnarray*}
\tilde R^{p,n}&=&\sum_{j=0}^{\ell-1} \tilde R^{p,n}_{\ell}(j) + o(1),
\qquad \tilde R^{p,n}_{\ell}(j) \ := \  \frac{1}{n-1}\sum_{k=0}^{[(n-2-j)/\ell]} \tilde \eta_n(k \ell +j).
\end{eqnarray*}
Then
$$
\E \left(\tilde R^{p,n}\right)^4 \  \le \ \ell^4 \max_{0\le j <\ell} \E
\left(\tilde R^{p,n}_{\ell}(j) \right)^4.  \label{R_bdd1}
$$
Write $\eta_n(k) $  as a (bounded) function in
standardized Gaussian variables:
\begin{equation}
\eta_n(k) \ = \  f_{k,n}\left({\mbf Y}_{\! n}(k)\right),
\end{equation}
where ${\mbf Y}_{\! n}(k) = (Y^{(1)}_n(k), Y^{(2)}_n(k)) \in
{\R}^2$,
\begin{eqnarray}
Y^{(1)}_{n}(k)&:=&\frac{\Delta^{p,n}_{k} X}{\sigma_{p,n}(k)},  \label{Y_11}\\
Y^{(2)}_{n}(k) &:=&-\frac{\Delta^{p,n}_{k} X}{\sigma_{p,n}(k)}
\frac{ \rho_{p,n}(k) }{ \sqrt{1 - \rho^2_{p,n}(k) }} +
\frac{\Delta^{p,n}_{k+1} X}{\sigma_{p,n}(k+1)} \frac{1}{\sqrt{1
- \rho^2_{p,n}(k)}}, \label{Y_21} \\
\mbox{and}~~~f_{k,n}\left(x^{(1)}, x^{(2)}\right) & :=&\psi\left(x^{(1)},
\frac{\sigma_{p,n}(k+1)}{\sigma_{p,n}(k)}\left(\rho_{p,n}(k) x^{(1)} +
\sqrt{1- \rho^2_{p,n}(k)}x^{(2)}\right)\right),  \label{f_kn1}
\end{eqnarray}
where $\sigma^2_{p,n}(k), \  \rho_{p,n}(k) $ are defined in
(\ref{rho_nk1}).  Then, for each $k$, ${\mbf Y}_{\! n}(k) := (Y^{(1)}_n(k), Y^{(2)}_n(k))$ has a standard Gaussian distribution in $\R^2 $ and
$\tilde \eta_n(k) \ = \  f_{k,n}\left({\mbf Y}_{\! n}(k)\right)- \E f_{k,n}\left({\mbf Y}_{\! n}(k)\right)$. Moreover, the vector $({\mbf Y}_{\! n}
(k \ell +j), \, k=0,1,\cdots, [(n-2-j)/\ell])\in \R^{2([(n-2-j)/\ell]+1)} $  is $\varepsilon-$standard provided $\ell $ is large enough.
Now Lemma \ref{lemgauss} can be used and it implies the bound (\ref{4mom1}) using Assumptions {\bf (A.1)}
and {\bf (A.2)$_p$}.
The details of this proof can be found in Bardet et Surgailis (2009).  \hfill
$\Box$

\subsection*{\it Sketch of the proof of Theorem \ref{TLCgauss}}

The proof of Theorem \ref{psgauss} uses the following central limit theorem for Gaussian subordinated multidimensional triangular arrays.
Theorem~\ref{tlcgauss} is proved in Bardet and Surgailis (2009). It
extends the earlier results in Breuer and Major~(1983) and Arcones~(1994). Below, similarly as in Lemma 1, ${\mbf X} \in \R^\nu$ designates
a standard Gaussian vector.

\begin{thm}\label{tlcgauss}
Let $\left({\mbf
Y}_{\! n}(k) \right)_{1\le k \le n, n \in {\N}} $ be a  triangular
array of standardized Gaussian vectors with values in ${\R}^\nu, \
{\mbf Y}_{\! n}(k) = (Y^{(1)}_n(k), \cdots, Y^{(\nu)}_n(k)), \  \E
Y^{(p)}_n(k) = 0, \  \E Y^{(p)}_n(k) Y^{(q)}_n(k) = \delta_{pq}$.
For a given integer $m\ge 1$, introduce the following assumption:
there exists a function $\rho:\N\to \R$ such that for any $1\le p, q \le \nu$,
\begin{eqnarray}
\forall (j, k)\in \{1,\cdots,n\}^2,~~~~\left|
\E Y^{(p)}_n(j) Y^{(q)}_n(k) \right|  &\le& |\rho(j-k)|~~~~\mbox{with} ~~\sum_{j\in {\Z}}
|\rho(j)|^m < \infty. \label{rho_dom}
\end{eqnarray}
Moreover, assume that for any $\tau \in [0,1]$ and
any $J \in \N^*$,
\begin{equation}
\left({\mbf Y}_{\! n}([n\tau]+j)\right)_{-J \le j \le J} \
\limiteloin \  \left({\mbf W}_{\! \tau}(j)\right)_{-J \le j \le J}
,   \label{Y_W}
\end{equation}
where $ \left({\mbf W}_{\! \tau}(j)\right)_{j\in {\Z}}$ is a
stationary Gaussian process taking values in ${\R}^\nu$ and
depending on parameter $\tau \in (0,1)$.
Let $\tilde f_{k,n} \in  \L^2_0({\mbf X}) \ (n \ge 1, 1 \le k \le n)$ be a
triangular array of functions all having Hermite rank at least $m$. Assume that there exists a $ \L^2_0({\mbf X})-$valued
continuous function $\tilde \phi_\tau, \tau \in [0,1]$, such that
\begin{equation}
\sup_{\tau \in [0,1]} \|\tilde f_{[\tau n],n} - \tilde \phi_{\tau}\|^2 =
\sup_{\tau \in [0,1]} \E (\tilde f_{[\tau n],n}({\mbf X}) -
\tilde \phi_{\tau}({\mbf X}))^2 \limiten 0. \label{f_n}
\end{equation}
Then, with $\displaystyle
\sigma^2=\  \int_0^1 \d \tau  \Big(\sum_{j\in {\Z}} \E \big [
\tilde \phi_\tau \left({\mbf W}_{\! \tau}(0)\right) \tilde \phi_\tau \left({\mbf
W}_{\! \tau}(j)\right)\big ]\Big) <\infty, $
\begin{equation}
n^{-1/2} \sum_{k=1}^n \tilde f_{k,n}\left({\mbf Y}_{\! n}(k)\right)
\limiteloin {\mathcal N}(0, \sigma^2). \label{CLT}
\end{equation}
\end{thm}
\noindent {\it Sketch of the proof of Theorem \ref{TLCgauss}.} It suffices to show that
\begin{equation} \label{Eproof}
\sqrt{n} \Big|\E R^{n,p} - \int_0^1 \Lambda_p (H(t)) \d t \Big|\limiten 0
\end{equation}
and
\begin{equation} \label{TLCproof}
\sqrt{n}\big (R^{p,n}  - \E R^{p,n})
\limiteloin {\mathcal N} \big (0,\int_0^1 \Sigma_p(H(\tau))\d \tau
 \big ).
 \end{equation}
The proof of (\ref{Eproof}) uses Assumption {\bf (A.1)'} or  (\ref{R_j11})  and the
easy fact that for Gaussian vectors $(Z^{(1)}_n,Z^{(2)}_n)\in \R^2, \, n\in \N $ with zero mean
$\E Z^{(i)}_n \equiv 0, \, i=1,2, n \in \N$ and $\E (Z^{(1)}_0)^2 =  \E (Z^{(2)}_0)^2 =1, \, |\E Z^{(1)}_0 Z^{(2)}_0| < 1 $
\begin{equation} \label{Zdist}
\big|\E \psi (Z^{(1)}_0 Z^{(2)}_0) - \E \psi(Z^{(1)}_n, Z^{(2)}_n) \big| \ \le \ C \sum_{i,j=1}^2
\big|\E Z^{(i)}_0 Z^{(j)}_0 - Z^{(i)}_n Z^{(j)}_n \big|.
\end{equation}
The proof of (\ref{TLCproof}) is deduced from Theorem \ref{tlcgauss} with the sequence of standardized Gaussian
vectors ${\mbf Y}_{\! n}(k) = (Y^{(1)}_n(k), Y^{(2)}_n(k)) \ (\nu
=2)$ given in (\ref{Y_11})-(\ref{Y_21}) and the centered functions
$$
\tilde f_{k,n}(x^{(1)},x^{(2)}):=f_{k,n}(x^{(1)},x^{(2)})-\E
f_{k,n}({\mbf Y}_{\! n}(k)),  \qquad \tilde \phi_\tau (x^{(1)},x^{(2)}):=  \phi_{\tau}(x^{(1)},x^{(2)})-\E
\phi_{\tau}({\mbf X})
$$
with $f_{k,n}:\R^2 \to \R$ given in (\ref{f_kn1}) and the (limit) function
$$
\phi_\tau (x^{(1)},x^{(2)}) :=  \psi\Big(x^{(1)},  \rho_p(H(\tau)) x^{(1)} +  \sqrt{1 - \rho^2_p(H(\tau))} x^{(2)}\Big).
$$
Thanks to symmetry properties of these functions, it is clear that the Hermite rank of
$\tilde f_{k,n}$ (for any $k$ and $n$) and $\tilde \phi_\tau $ (for any $\tau \in
[0,1]$) is $m=2$. Using Assumptions {\bf (A.1)'}
and {\bf (A.2)$_p$} (with $\theta = 0, \gamma> 1/2$), one can show that the conditions of Theorem \ref{tlcgauss} are satisfied for  the above
$\tilde
f_{n,k}, \tilde \phi_\tau $ and the limit process $({\mbf W}_\tau
(j))_{j\in {\Z}} $ in (\ref{Y_W}) is written in terms of increments of fBm
$(B_{H(\tau)}(j))_{j \in \Z}$:
$$
{\mbf W}_\tau(j) \ := \  \left( \Delta^p_1B_{H(\tau)}(j), \left(-
\rho_p(H(\tau)) \Delta^p_1 B_{H(\tau)}(j) + \Delta^p_1
B_{H(\tau)}(j+1)\right)/\sqrt{1- \rho_p^2(H(\tau))} \right)
$$
having standardized uncorrelated components. The details of this proof can be found in Bardet and Surgailis (2009).  \hfill
$\Box$

\subsection*{\it Proof of Proposition \ref{trends}}

\noindent (i)
Denote $\tilde Z_t := \alpha(t) X_t, \  R^{n,p} :=  R^{n,p}(Z),  \tilde R^{n,p} := R^{n,p} (\tilde Z)$. Clearly, part (i) (the CLT in (\ref{TLC}) for trended process $Z = \alpha X + \beta $) follows from
the following relations:
\begin{eqnarray}
&&\sqrt{n}\, \big (\tilde R^{p,n}  - \E \tilde R^{p,n})
\ \limiteloin \ {\mathcal N}\big (0,\int_0^1 \Sigma_p(H(\tau))\d \tau
 \big ),\label{TLC1} \\
&&\sqrt{n} \, \Big|\E \tilde R^{p,n}- \int_0^1 \Lambda_p(H(t)) \d t\Big|\ \limiten \ 0, \label{TLC2} \\
&& \sqrt{n}\,  \big|\E R^{n,p} - \E \tilde R^{n,p}\big| \ \limiten \ 0, \label{TLC3} \\
&&n \, \Var(R^{n,p} - \tilde R^{n,p})\ \limiten \ 0.  \label{TLC4}
\end{eqnarray}
The central limit theorem in (\ref{TLC1}) follows from Theorem \ref{TLCgauss} and Remark \ref{weakTLC} since Assumption {\bf (A.2)$_p$}
and the convergences (\ref{sigma1})-(\ref{corr2}) (with $c(t)$ replaced by $\alpha^2(t) c(t)$) can be easily
verified for the process $\tilde Z = \alpha X $.

Let us turn to the proof of (\ref{TLC1})-(\ref{TLC4}). For concreteness, let $p=2$ in the rest of the proof. Since
$\Delta^{2,n}_k \tilde Z = \alpha(k/n) \Delta^{2,n}_k X + 2 \alpha'(k/n) n^{-1} \Delta^{1,n}_{k+1} X + O(1/n^2) X((k+2)/n) $, it
follows easily from  Assumption ${\bf (A.1)'}$ and (\ref{R_j11}) (for $X$) that
$$
\E \big[\Delta^{2,n}_{[nt]} \tilde Z \, \Delta^{2,n}_{j+[nt]} \tilde Z\big] \ = \  \alpha^2(t) \E
\big[\Delta^{2,n}_{[nt]} X \, \Delta^{2,n}_{j+[nt]} X \big]
+ O(n^{-2H(t)-1} + n^{-H(t) -2})
$$
implying (\ref{R_j11}) for $\tilde Z $ (with $c(t)$ replaced by $\alpha^2(t) c(t)$). Whence and using (\ref{R_j11}) and (\ref{Zdist}),  relation (\ref{TLC2}) follows
similarly as (\ref{Eproof}) above.

The proofs of (\ref{TLC3})-(\ref{TLC4}) uses the following bounds from Bru\v zait\.e and Vai\v ciulis (2008, Lemma 1). Let $(U_1, U_2)\in \R^2 $
be a Gaussian vector with zero mean, unit variances and a correlation coefficient $\rho, \ |\rho| < 1 $. Then for any  $b_1, b_2 \in \R$
and any $1/2 < a_1, a_2 < 2 $
 \begin{eqnarray}
&&\big|\E \psi(a_1 U_1 + b_1, a_2 U_2 + b_2) - \E \psi(a_1 U_1, a_2 U_2)\big|\ \le\  C(b_1^2 + b_2^2), \label{U1}\\
&&\big|\E U_i \psi(a_1 U_1 + b_1, a_2 U_2 + b_2)\big|\ \le \ C(|b_1| + |b_2|), \quad i=1,2, \label{U2}
\end{eqnarray}
where the constant $C$ depends only on $\rho $ and does not depend on $a_1, a_2, b_1, b_2$. Using (\ref{U1})
, Assumption ${\bf (A.1)'}$ and
the fact that  $|\Delta^{2,n}_k \beta| \le C n^{-2}$ we obtain
\begin{equation} \label{ERconv}
 \big|\E R^{n,2} - \E \tilde R^{n,2}\big| \ \le \  C  n^{-1} \sum_{k=0}^{n-3} \big(n^{H(k/n)} n^{-2} \big)^2 \ = \ O(n^{-2}),
\end{equation}
proving (\ref{TLC3}).

To prove (\ref{TLC4}), write $R^{n,2} = \frac{1}{n-2} \sum_{k=0}^{n-3} \eta_n(k), \
\tilde R^{n,2} = \frac{1}{n-2} \sum_{k=0}^{n-3} \tilde \eta_n(k)$, where
\begin{eqnarray*}
\eta_n(k)&:=&\psi\big(\Delta^{2,n}_k Z, \Delta^{2,n}_{k+1} Z\big)\ = \  f_{n,k} \big(\tilde Y^{(1)}_n(k) + \mu^{(1)}_n(k), \tilde Y^{(2)}_n(k) + \mu^{(2)}_n(k)\big),  \\
\tilde \eta_n(k)&:=&\psi\big(\Delta^{2,n}_k \tilde Z, \Delta^{2,n}_{k+1} \tilde Z\big) \ = \ f_{n,k} \big(\tilde Y^{(1)}_n(k), \tilde Y^{(2)}_n(k)\big),
\end{eqnarray*}
where standardized increments $\tilde Y^{(i)}_n(k), i=1,2 $ are defined as in (\ref{Y_11})-(\ref{Y_21}) with $X$ replaced by $\tilde Z$,
$f_{n,k} $ are defined in (\ref{f_kn1}), and
$$
\mu^{(1)}_n(k) := \frac{\Delta^{2,n}_k \beta}{\sigma_{2,n}(k)}, \quad
\mu^{(2)}_n(k) := - \frac{\Delta^{2,n}_{k} \beta}{\sigma_{2,n}(k)}\frac{\rho_{2,n}(k)}{\sqrt{1- \rho^2_{2,n}(k)}}
+ \frac{\Delta^{2,n}_{k+1} \beta}{\sigma_{2,n}(k+1)}\frac{1}{\sqrt{1- \rho^2_{2,n}(k)}}.
$$
Note, the $\tilde \eta_n(k)$'s and $\eta_n(k)$'s have the Hermite rank $\ge 2 $ and $\ge 1$, respectively, since the Hermite coefficients of order $1$
$c^{(i)}_n(k) := \E \big[Y^{(i)}_n(k) \eta_n(k)\big], \, i=1,2$ of the $\eta_n(k)$'s are not  zero in general. Using the bound in (\ref{U2}) and 
$|\Delta^{2,n}_k \beta| \le C n^{-2}$ we obtain
\begin{equation} \label{cbdd}
|c^{(i)}_n(k)| \le C\big(|\mu^{(1)}_n(k)| + |\mu^{(1)}_n(k)|\big)  \le Cn^{H(k/n)-2}.
\end{equation}
Split $\eta_n(k)-
\tilde \eta_n(k) - (\E \eta_n(k) - \E \tilde \eta_n(k)) = \zeta'_n(k) + \zeta''_n(k) $, where
\begin{equation} \label{zetadef}
\zeta'_n(k):=\sum_{i=1}^2 c^{(i)}_n(k) Y^{(i)}_n(k), \quad
\zeta''_n(k) := \eta_n(k)- \zeta'_n(k) -
\tilde \eta_n(k) - (\E \eta_n(k) - \E \tilde \eta_n(k)).
\end{equation}
Then $\Var(R^{n,2}- \tilde R^{n,2}) \le 2(J'_n + J''_n)$,  where
$J'_n := \E \Big(\frac{1}{n-2} \sum_{k=0}^{n-3} \zeta'_n(k)\Big)^2, \
J''_n:= \E \Big(\frac{1}{n-2} \sum_{k=0}^{n-3} \zeta''_n(k)\Big)^2 $. From
(\ref{cbdd}) and Cauchy-Schwartz inequality, it follows $J'_n \le \frac{C}{(n-2)^2}\big ( (n-2)\sum_{k=0}^{n-3} n^{2H(k/n)-4} \big) = O(1/n^2) = o(1/n)$. Since the
$\zeta''_n(k)$'s have Hermite rank $\ge 2$, we can use the
argument in the proof of Theorem \ref{TLCgauss} together with Arcones' inequality (Arcones, 1994) and the easy fact
that $\E (\zeta''_n(k))^2 \to 0 \ (n \to \infty, \ k=0, \cdots, n-3)$, to conclude that $J''_n = o(1/n) $. This proves
(\ref{TLC4}) and completes the proof of part (i).

\smallskip

\noindent (ii) Similarly as in the proof of part (i), we shall restrict ourselves to the case $p=2$ for concreteness. We use
the same notation and the decomposition of $R^{2,n}$ as in part (i):
\begin{eqnarray} \label{Rdecomp}
R^{2,n} - \int_0^1 \Lambda_p(H(\tau))\d \tau&=&\Big(\tilde R^{2,n} - \int_0^1 \Lambda_p(H(\tau))\d \tau\Big) + \big(\E R^{n,p} - \E \tilde R^{n,p}\big) + Q_n,
 \end{eqnarray}
where
\begin{eqnarray*}
Q_n&:=&R^{2,n} - \tilde R^{2,n} - \big(\E R^{n,p} - \E \tilde R^{n,p}\big)\ = \ \frac{1}{n-2}\sum_{k=0}^{n-3} \zeta'_n(k) + \frac{1}{n-2}\sum_{k=0}^{n-3} \zeta''_n(k) \ =: \ Q'_n + Q''_n;
 \end{eqnarray*}
see (\ref{zetadef}). Here, the a.s. convergence to zero of the first term on the r.h.s. of (\ref{Rdecomp}) follows from
Theorem \ref{psgauss} since the conditions of this theorem for $\tilde Z$ are easily verified. The convergence to zero
of the second term on the r.h.s. of (\ref{Rdecomp}) follows similarly as in (\ref{ERconv}), with the difference that
$|\Delta^{2,n}_k \beta| \le C n^{-1}$ since $\beta \in {\cal C}^1 [0,1]$ and therefore
\begin{equation} \label{EER}
\big|\E R^{n,2} - \E \tilde R^{n,2}\big| \ \le \  C  n^{-1} \sum_{k=0}^{n-3} \big(n^{H(k/n)} n^{-1} \big)^2 \ = \
O(n^{-2(1- \sup_{t\in [0,1]} H(t))})
= o(1).
\end{equation}
The proof of $Q''_n \limitepsn 0$ mimics that of Theorem \ref{psgauss} and relies on the fact that the
$\zeta''_n(k)$'s have the Hermite rank $\ge 2$ (see above). Relation
$Q'_n \limitepsn 0$ follows by the gaussianity of $Q'_n$ and  $\E (Q'_n)^2 = O(n^{-\delta}) $ for some $\delta > 0$.
Since $|c^{(i)}_n(k)| \le   Cn^{H(k/n)-1} $ in view of the first inequality in (\ref{cbdd}) and
$|\Delta^{2,n}_k \beta| \le C n^{-1}$, then $\E (Q'_n)^2 \le Cn^{-2(1- \sup_{t\in [0,1]} H(t))} $ similarly
as in (\ref{EER}) above. This proves part (ii) and the proposition. \hfill
$\Box$

\subsection*{\it Proof of Corollary \ref{corr0}}

(a) The argument at the end of the proof of Property 4.1 shows that
$V$ satisfies Assumption ${\bf (A.1)}'$, while ${\bf (A.2)}_p$
follows from Property \ref{prop_V}. Then the central limit theorem
in (\ref{TLCstat}) follows from Theorem 4.2.

\smallskip

\noindent (b) In this case, ${\bf (A.2)}_p$ follows from Property
\ref{prop_f}. Instead of verifying ${\bf (A.1)}'$, it is simpler to
directly  verify condition (\ref{R_j11}) which suffices for the
validity of the statement of Theorem 4.2. Using  $f(\xi) = c
\xi^{-2H-1}\big (1+o(\xi^{-1/2})\big ) \ (\xi\to \infty)$, similarly
as in the proof of Property \ref{prop_f}  for $j \in \N^*$ one
obtains
\begin{eqnarray*}
&&\sqrt n \, \Big |n^{-2H} \E \big [\Delta^{p,n}_0 X \,
\Delta^{p,n}_{j} X \big ] - c \, 2^{1+p}\int_{0}^\infty (1- \cos
(x))^p \cos(xj) \, x^{-2H-1} \d x \Big | \\ &&= \  2^{1+p}\Big|
\int_{0}^\infty (1- \cos (x))^p \cos(xj) \times  \sqrt n \, \big (
n^{2H+1}(f(nx)- c(nx)^{-2H-1})\big) \d x\Big| \   \limiten 0
\end{eqnarray*}
by the Lebesgue dominated convergence theorem,  since $\int_{0}^\infty \big |(1- \cos (x))^p
\cos(xj) \, x^{-2H-3/2}\big | \d x <\infty$. Therefore condition
(\ref{R_j11}) is satisfied and Theorem
\ref{TLCgauss} can be applied. \\
Finally, the function $H \mapsto
\Lambda_2(H)$ is a ${\cal C}^1(0,1)$ bijective function and from
the Delta-method (see for instance Van der Vaart, 1998), the
central limit theorem in (\ref{TLCstat2}) is shown. \hfill $\Box$

\subsection*{\it Proof of Lemma \ref{lemmaFG}}

(i) We use the following inequality (see, e.g. Ibragimov and Linnik, 1971, Theorem 1.5.2):
\begin{equation}\label{E}
\|F_n-G_\alpha\|_\infty \le \frac{1}{2\pi} \int_{\R} \Big|\frac{f_n(\theta)-g_\alpha(\theta)}{\theta}\Big| {\d}\theta,
\end{equation}
where $f_n, g_\alpha $ are the characteristic functions of $F_n, G_\alpha$, respectively.
According to (\ref{K}) and the definition of $K_1$,
$$
f_n(\theta) - g_\alpha(\theta)\ = \ g_\alpha (\theta) \big(\e^{-2n^{-1} I_n(\theta) } -1 \big),
$$
where, with $v := \theta n^{1/\alpha}$,
$$
I_n(\theta)\ := \ \int_0^\infty
{\rm Re}(\e^{{\i} u \,v } -1) {\d} K_1(u) \ = \ \int_0^{1/v} \cdots  + \int_{1/v}^\infty  \cdots
\ =: \ I_1 + I_2.
$$
If $v > 1 $ then integrating by parts and using (\ref{KK}),
\begin{eqnarray*}
|I_2|&\le&2|K_1|(1/v) \ = \  O(v^{(\alpha - \beta)_+}),  \\
|I_1|&=&\Big|K_1(1/v){\rm Re}(\e^{{\i}} -1) - \int_0^{1/v} K_1(u)\, v \, \sin( u v)\, {\d}u  \Big|\\
&\le&2|K_1(1/v)| + v^2 \int_0^{1/v} u |K_1(u)|{\d}u \\
&=&O(v^{(\alpha - \beta)_+}) + v^2 (1/v)^{2-(\alpha - \beta)_+}  =O(v^{(\alpha - \beta)_+}).
\end{eqnarray*}
Next, if $v \le 1 $, then similarly as above
\begin{eqnarray*}
|I_2|&\le&2|K_1|(1/v) \ = \  O(v^{\delta}),  \\
|I_1|&=&\Big|K_1(1/v){\rm Re}(\e^{{\i}} -1) - \int_0^{1/v} K_1(u) \, v \, \sin( u v)\, {\d}u \Big|\\
&\le&2|K_1(1/v)| + v^2 \int_0^{1/v} u |K_1(u)|{\d}u \\
&=&O(v^\delta) + Cv^2 \Big(\int_0^1 u u^{-(\alpha - \beta)_+} {\d}u + \int_1^{1/v}
u^{1-\delta} {\d}u \Big)=O(v^{\delta \wedge 2}).
\end{eqnarray*}
Therefore, for some constant $C$,
$$
|I_n(\theta)| \le  C\cases{(\theta n^{1/\alpha})^{(\alpha - \beta)_+}, &$\theta > 1/n^{1/\alpha}$, \cr
(\theta n^{1/\alpha})^{\delta \wedge 2}, &$\theta \le 1/n^{1/\alpha}$.\cr}
$$
Moreover, $g_\alpha(\theta) = \e^{-c_1|\theta|^\alpha}$ for some $c_1 >0$. Using these facts
and (\ref{E}) we obtain
\begin{eqnarray*}
\|F_n-G_\alpha\|_\infty&\le&C\Big(\int_0^{1/n^{1/\alpha}} n^{-1} (\theta n^{1/\alpha})^{\delta \wedge 2} \frac{{\d}\theta}{\theta} +
+ \int_{1/n^{1/\alpha}}^\infty n^{-1} (\theta n^{1/\alpha})^{(\alpha - \beta)_+} \e^{-c_1 |\theta|^\alpha}\frac{{\d}\theta}{\theta}\Big)\\
&\le&Cn^{-1} + C \frac{1}{n^{1 - (\alpha - \beta)_+/\alpha}} \Big(1 + \int_{1/n^{1/\alpha}}^1 \theta^{(\alpha - \beta)_+ -1} {\d}\theta \Big).
\end{eqnarray*}
The above bound easily yields (\ref{FG}).

\smallskip

\noindent (ii) Follows similarly using (\ref{E}) and the argument in (i) with $K_1$ replaced by $K$. \hfill $\Box$

\vskip.4cm

\noindent{\large \it Computation of $\lambda(r)$}

\vskip.3cm

\noindent From the definition
of $\lambda(r)$ and the change of variables $x_1=a \cos\phi$,
$x_2=a \sin \phi$, with $|r|<1$,
\begin{eqnarray*}
\lambda(r)&=&\frac 1 {2\pi \sqrt{1-r^2}} \int _{\R^2} \frac
{|x_1+x_2|}{|x_1|+|x_2|} \e^{-\frac 1 {2(1-r^2)}
(x_1^2-2rx_1x_2+x_2^2)} \d x_1\d x_2 \\
&=& \frac {\sqrt{1-r^2}} \pi \int_0^\pi \frac
{|\cos\phi+\sin\phi|}{(|\cos \phi| + |\sin \phi|)
(1-r\sin(2\phi))} \d\phi \\
&=:&I_1+I_2,
\end{eqnarray*}
where
\begin{eqnarray*}
I_1&=&\frac {\sqrt{1-r^2}} \pi \int_0^{\pi/2} \frac 1 {
1-r\sin(2\phi)} \d\phi\\
&=&\frac {\sqrt{1-r^2}} \pi \int_0^{\infty} \frac 1 { 1+t^2-2rt}
\d t\ = \ \frac 1 2 +\frac 1 \pi \arctan \Big ( \frac r {\sqrt{1-r^2}} \Big)\ =\  \frac 1 \pi \arccos(-r); \\
I_2&=&\frac {2\sqrt{1-r^2}} \pi \int_0^{\pi/4} \frac {\cos
\phi-\sin \phi} { (\cos \phi + \sin \phi)(1+r\sin(2\phi))}
\d\phi\ =\  \frac {\sqrt{1-r^2}} \pi \int _0 ^1 \frac
{1-t}{(1+t)(1+2rt+t^2)}\d t \\
&=&\frac {\sqrt{1-r^2}} {\pi(1-r)}\log \Big ( \frac 2 {r+1}\Big ).
\end{eqnarray*}
The function
$\lambda(r) $ is monotone increasing on $[-1,1]; \ \lambda(1)
=1, \ \lambda(-1) = 0.$ It is easy to check that
\begin{equation}
\rho_{1}(H) \ = \ 2^{2H-1}-1, \qquad \rho_{2}(H) =
\frac{-3^{2H}+2^{2H+2}-7}{8-2^{2H+1}}   \label{rho_1}
\end{equation}
are monotone increasing functions; $ \rho_{1}(1) = 1,
\rho_{2}(1) =0 $ so that $\Lambda_{p}(H) = \lambda
(\rho_{p}(H)) $ for $p=1,2 $ is also monotone for $H \in (0,1)$.

\vskip.4cm

\noindent {\large \it Expression and graph of $s_2(H)$}

\vskip.3cm

\noindent From the Delta-method, $\displaystyle{s_2^2(H)=\Big [ \frac
{\partial}{\partial x}  (\Lambda_2)^{-1}(\Lambda_2(H))\Big
]^2\Sigma_2(H) }$ and therefore $$ \displaystyle{s_2^2(H)=\Big (
\frac {\pi \, (8-2^{2H+1})^2 (1-\rho_2(H))
\sqrt{1-\rho_2^2(H)}}{\big (\log 2 -\log (1+\rho_2(H))\big ) \big
(2^{2H+2}9\log 2 -3^{2H}16\log 3 +6^{2H} 4 \log (3/2)\big )}\Big )^2
\, \Sigma_2(H)},
$$
with the approximated graph (using the numerical values of
$\Sigma_2(H)$ in Stoncelis and Vai\v ciulis (2008)):
\begin{center}
\begin{figure}[h]
\begin{center}
\includegraphics[width=10 cm,height=4cm]{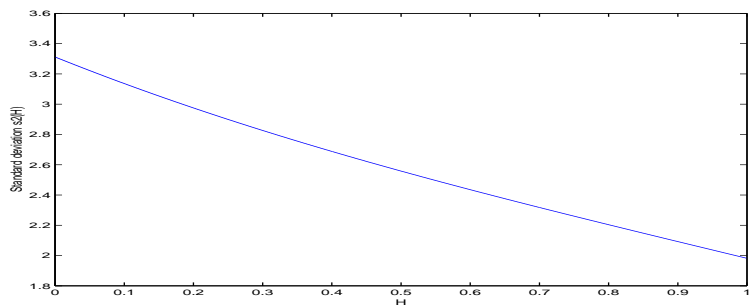}
\caption{The graph of $\sqrt{s_2^2(H)}$ } \label{Figure4}
\end{center}
\end{figure}
\end{center}

\noindent {\bf Acknowledgements.} The authors are grateful to the referee for many relevant
suggestions and comments which helped to improve the contents of the paper.  The authors also want to thank Antoine
Ayache
for a discussion, and Mindaugas Stoncelis and Marijus Vai\v ciulis
for providing us with the numerical computations
for the graphs in Figures \ref{Figure3}, \ref{Figure2}  and \ref{Figure4}.


\begin{thebibliography}{9}



\bibitem{jacod} Ait Sahalia, Y. and Jacod, J. (2009) Estimating the degree of activity of jumps in high
frequency financial data. {\em Ann. Statist.} {\bf 37}, 2202--2244.


\bibitem{arcones} Arcones, M.A. (1994) Limit theorems for nonlinear
functionals of a
    stationary Gaussian sequence of vectors.
{\em Ann. Probab.} {\bf 22}, 2242--2274.


\bibitem{Ayache2} Ayache, A., Benassi, A., Cohen, S. and L\'evy V\'ehel,
J. (2005) Regularity and identification of generalized
multifractional Gaussian processes. In: Lecture Notes in Math. vol.
1857, pp. 290--312, Springer, Berlin.



\bibitem{Bardet} Bardet, J.-M. and Bertrand, P. (2007) Definition, properties
and wavelet analysis of multiscale fractional Brownian motion,
{\em Fractals} {\bf 15}, 73--87.

\bibitem{ba-sur1} Bardet, J.-M. and Surgailis, D. (2009)  A central limit theorem for triangular arrays of nonlinear functionals of Gaussian vectors. Preprint.

\bibitem{ba-sur2} Bardet, J.-M. and Surgailis, D. (2010) Nonparametric estimation of the local Hurst function of multifractional processes. Preprint.

\bibitem{bel} Belomestny, D. (2010) Spectral estimation of the fractional order of a Lévy process. {\it Ann. Statist.} {\bf 38}, 317--351.

\bibitem{Benassi3} Benassi, A., Cohen, S. and Istas, J. (1998) Identifying the
multifractal function of a Gaussian process. {\em Statist. Probab.
Lett.} {\bf 39}, 337--345.



\bibitem{Benassi} Benassi, A., Jaffard, S. and Roux, D. (1997) Gaussian
processes and pseudodifferential elliptic operators. {\em Rev. Mat. Iberoamericana}   {\bf 13}, 19--89.




\item{Berk, K.N.} (1973) A central limit theorem for
$m-$dependent random variables with unbounded $m$.  {\em Ann. Probab.} {\bf 1}, 352--354.


\bibitem{Breuer} Breuer, P. and Major, P. (1983) Central limit theorems
for nonlinear functionals of Gaussian fields. {\em J.  Multiv.
Anal.} {\bf 13}, 425--441.


\bibitem{Bruzaite} Bru\v zait\.e, K. and Vai\v ciulis, M. (2008) The increment ratio statistic
under deterministic trends. {\em Lithuanian  Math. J.} {\bf 48}, 256--269.

\bibitem{Coeurjolly1} Coeurjolly, J.-F. (2001) Estimating the parameters of a fractional
Brownian motion by discrete variations of its sample paths.  {\em Stat. Inference Stoch. Process.} {\bf 4}, 199--227.

\bibitem{Coeurjolly2} Coeurjolly, J.-F. (2005) Identification of multifractional Brownian motion.
{\em Bernoulli} {\bf 11}, 987--1008.

\bibitem{Coeurjolly3} Coeurjolly, J.-F. (2007) Hurst exponent estimation of locally self-similar Gaussian processes using sample quantiles.
{\em  Ann. Statist.}  {\bf  36}, 1404--1434.



\bibitem{Cramer} Cram\'er, H. and Leadbetter, M.R.
(1967) {\em Stationary and Related Stochastic Processes. Sample
Function Properties and Their Applications.} Wiley, New York.

\bibitem{csorgo} Cs\"org\H o, M. and Mielniczuk, J. (1996) The empirical process of a short-range dependent stationary sequence
under Gaussian subordination. {\em Probab. Th. Rel. Fields} {\bf 104}, 15--25.


\bibitem{Dahlhaus} Dahlhaus R. (1989) Efficient parameter estimation
for self-similar processes. {\em Ann. Statist.} {\bf 17},
1749--1766.

\bibitem{Dobrushin} Dobrushin, R.L. (1980) Automodel generalized random fields and their
renorm group. In: Dobrushin, R.L. and Sinai, Ya.G. (eds.){\em
Multicomponent Random Systems}, pp. 153--198. Dekker, New York.




\bibitem{Falconer1} Falconer, K. (2002) Tangent fields and the local
structure of random fields. {\em J. Theor. Probab.} {\bf 15},
731--750.

\bibitem{Falconer2} Falconer, K. (2003) The local structure of
random processes. {\em J. London Math. Soc.} {\bf 67}, 657--672.

\bibitem{Feuerverger} Feuerverger, A., Hall, P. and Wood, A.T.A. (1994)
Estimation of fractal index and fractal dimension of a Gaussian
process by counting the number of level crossings. {\em J. Time
Series Anal.} {\bf 15}, 587--606.

\bibitem{Fox} Fox, R. and Taqqu, M.S. (1986) Large-sample
properties of parameter estimates for strongly dependent Gaussian
time series. {\em Ann. Statist.} {\bf 14}, 517--532.

\bibitem{Gikhman} Gikhman, I.I. and Skorohod, A.V. (1969) {\em Introduction to the
Theory of Random Processes.} Saunders, Philadelphia.

\bibitem{Gugu} Gugushvili, S. (2008) Nonparametric estimation of the
characteristic triplet of a discretely observed Lévy process. {\em J. Nonparametr. Stat.} {\bf 21}, 321--343.

\bibitem{Guyon} Guyon, X. and Leon, J. (1989) Convergence en loi
des H-variations d'un processus gaussien stationnaire. {\em Ann.
Inst. Poincar\'e} {\bf 25}, 265--282.

\bibitem{Hall} Hall, P. and Wood, A. (1993) On the performance of box-counting estimators of fractal dimension.
{\em Biometrika}  {\bf 80}, 246--252.

\bibitem{HoSun} Ho, H.-C. and Sun, T.C. (1987)
A central limit theorem for non-instantaneous filters of a stationary Gaussian. {\em J. Multiv. Anal.} {\bf 22}, 144--155.


\bibitem{IbLi} Ibragimov, I.A. and Linnik. Y.V. (1971) {\em Independent and Stationary Sequences of
Random Variables.} Wolters-Noordhoff, Groningen.

\bibitem{Istas} Istas, J. and Lang, G. (1997) Quadratic variations and
estimation of the local H\"older index of a Gaussian process. {\em
Ann. Inst. Poincar\'e} {\bf 33}, 407--436.




\bibitem{Neum}  Neumann, M. and Reiß, M. (2009) Nonparametric estimation for Lévy processes from
low-frequency observations. {\em Bernoulli} {\bf 15}, 223--248.

\bibitem{Pel} Peltier, R. and Lévy-Vehel, J. (1994) A new method for estimating the parameter of fractional Brownian motion.
Preprint INRIA. Preprint available on {\tt
http://hal.inria.fr/docs/00/07/42/79/PDF/RR-2396.pdf}.

\bibitem{Phil1} Philippe, A., Surgailis, D. and Viano,  M.-C. (2006) Invariance principle
for a class of non stationary processes with long memory. {\em C. R.
Acad. Sci. Paris Ser. 1 } {\bf 342}, 269--274.

\bibitem{Phil2} Philippe, A., Surgailis, D. and Viano,  M.-C. (2008)  Time-varying
fractionally integrated processes with nonstationary long memory.
{\em  Th. Probab. Appl.}   {\bf 52}, 651--673.


\bibitem{Sato} Sato, K.-I. (1999) {\em L\'evy Processes and Infinitely Divisible Distributions.} Cambridge Univ. Press, Cambridge.

\bibitem{Shil} Shilov,  G.E. and  Gurevich, B.L. (1967) {\em Integral, Measure and Derivative}
(in Russian). Nauka, Moscow.


\bibitem{Ston} Stoncelis, M. and Vai\v ciulis, M. (2008) Numerical
approximation of some infinite gaussian series and integrals. {\em Nonlinear Analysis: Modelling and Control} {\bf 13},  397--415.






\bibitem{Su} Surgailis, D.(2008) Nonhomogeneous fractional integration and
multifractional processes. {\em Stochastic Process. Appl.} {\bf
118}, 171--198.

\bibitem{SuTeVa} Surgailis, D., Teyssi\`ere, G. and Vai\v ciulis, M. (2008)
The increment ratio statistic. {\em J. Multiv. Anal.} {\bf 99}, 510--541.

\bibitem{Taqqu1} Taqqu, M.S. (1977)  Law of the iterated logarithm for sums of
non-linear functions of Gaussian variables that exhibit a long range
dependence.  {\em Z. Wahrsch. Verw. Geb.} {\bf 40}, 203--238.

\bibitem{Vai} Vai\v ciulis, M. (2009) An estimator of the tail index based on increment ratio statistics.
{\em Lithuanian Math. J.} {\bf 49}, 222--233.


\bibitem{Van} Van der Vaart, A. (1998) {\em Asymptotic Statistics}.
Cambridge Series in Statistical and Probabilistic Mathematics,
Cambridge.
\end{thebibliography}
\end{document}